\documentclass[a4paper,reqno]{amsart}
\usepackage{amsmath}
\usepackage{amssymb}
\usepackage{amsthm}
\usepackage{amsfonts}
\usepackage{mathrsfs}
\usepackage{graphicx}
\usepackage{booktabs}
\usepackage{cases}
\usepackage{indent}
\usepackage{subfigure}
\makeatletter
	
	\@addtoreset{equation}{section}
\makeatother
\makeatletter
	
	\@addtoreset{figure}{section}
\makeatother
\makeatletter
	
	\@addtoreset{table}{section}
\makeatother
\theoremstyle{definition}
\newtheorem{Def}{Definition}[section]
\theoremstyle{definition}
\newtheorem{Prob}{Problem}[section]
\theoremstyle{definition}
\newtheorem{Alg}{Algorithm}[section]
\theoremstyle{plain}
\newtheorem{Thm}{Theorem}[section]
\theoremstyle{plain}
\newtheorem{Prop}{Proposition}[section]
\theoremstyle{plain}
\newtheorem{Lem}{Lemma}[section]
\theoremstyle{remark}
\newtheorem{Rem}{Remark}[section]
\title[finite element approximation of the Stokes problem]{On a finite element approximation of the Stokes problem under leak or slip boundary conditions of friction type}
\author[T. Kashiwabara]{Takahito Kashiwabara}
\address{Graduate~School of Mathematical Sciences, The University of Tokyo, 3-8-1 Komaba, Meguro, Tokyo 153-8914, Japan}
\email{tkashiwa@ms.u-tokyo.ac.jp}
\keywords{Finite element method, Stokes equations, Boundary conditions of friction type, Variational inequality, Uzawa algorithm}
\subjclass[2010]{65N30, 35Q30, 35J87}

\begin{document}
\begin{abstract}
	A finite element approximation of the Stokes equations under a certain nonlinear boundary condition, namely, the slip or leak boundary condition of friction type, is considered.
	We propose an approximate problem formulated by a variational inequality, prove an existence and uniqueness result, present an error estimate, and discuss a numerical realization
	using an iterative Uzawa-type method. Several numerical examples are provided to support our theoretical results.
\end{abstract}

\maketitle

\section{Introduction}\label{Sec1}
We consider the motion of an incompressible fluid in a bounded two-dimensional domain with some nonlinear boundary conditions, specified as the
\emph{slip boundary condition of friction type} (SBCF) or the \emph{leak boundary condition of friction type} (LBCF).
These boundary conditions were introduced by H.~Fujita in \cite{F94}, and subsequently, many studies have focused on the properties of the solution, for example, existence, uniqueness, regularity, and continuous dependence on data, for the Stokes and Navier-Stokes equations under such boundary conditions.
Details can be referred to in \cite{F94} itself or in \cite{S04}, \cite{RT07}, \cite{S07}, and \cite{ALL09}, among others.
Similar types of nonlinear boundary condition, such as subdifferential boundary condition or Tresca boundary condition, have been reported in \cite{K00}, \cite{C03}, and \cite{BB03}, among others.

The frictional boundary conditions under consideration have been successfully applied to some flow phenomena in environmental and medical problems
such as oil flow over or beneath sand layers and blood flow in the thoracic aorta. Such applications have been discussed in \cite{KFS98}, \cite{SK04}, \cite{SUR09}, and \cite{SU10}.
In these works, the finite difference method is used for discretization, and theoretical considerations such as convergence are not addressed.

On the other hand, few studies have focused on the theoretical analysis of numerical methods for these boundary conditions, even if restricted to the Stokes problem.
For example, Li and Li \cite{LL08} proposed a finite element approximation combined with a penalty method for the Stokes equation with SBCF.
They proved the optimal order error estimate; however, they did not focus on a numerical realization of their finite element approximation.

The purpose of this work is to construct a comprehensive theory of the finite element method applied to flow problems with SBCF and LBCF,
including all of the existence and uniqueness result, error analysis, and numerical implementation.
In doing so, herein, we restrict our consideration to the stationary Stokes equation in a two-dimensional polygon.

The remainder of this paper is organized as follows. In Section 2, we review the results for the continuous problems described in \cite{F94}.
Weak formulations by an elliptic variational inequality for SBCF and LBCF are also presented.
In Section 3, we prepare the finite element framework using the so-called P2/P1 element, and state several technical lemmas.

Section 4 is devoted to the study of approximate problems for SBCF.
We propose the discretized variational inequality problem, proving the existence and uniqueness of a solution.
In the error analysis, we first derive a primitive result of the convergence rate $O(h^{\min\{\epsilon,1/4\}})$ under the $H^{1+\epsilon}$-$H^{\epsilon}$ regularity assumption with $0<\epsilon\le2$.
Second, we show that it is improved to $O(h^{\min\{\epsilon,1\}})$ under the additional hypothesis of good behavior of the sign of the tangential velocity component on the boundary where SBCF is imposed.
A sufficient condition to obtain $O(h^\epsilon)$, which is of optimal order when $\epsilon=2$, is also considered.
Finally, we propose an iterative Uzawa-type algorithm to perform numerical computations, and prove that the iterative solution indeed converges to the desired approximate solution.

Section 5 is devoted to the study of approximate problems for LBCF, in a manner similar to Section 4.
However, it should be noted that unlike in the case of SBCF, we have to explicitly deal with an additive constant for the pressure.
As a result, sometimes there exist multiple solutions for the pressure, especially its additive constant; other times it is uniquely determined.
Moreover, in an error analysis, we can only obtain the convergence rate $O(h^{\min\{\epsilon/2,1/4\}})$, because of the error of the additive constant of the pressure.
If we can eliminate the influence of this error, the same rate-of-convergence as in the case of SBCF is realized.

In Section 6, several numerical examples are provided to support our theory.
We observe that the results of our computation capture the features of SBCF and LBCF and that the numerically calculated errors decrease at $O(h^2)$ for both.
Section 7 presents the conclusions and discusses some future works.

The author learned about Ayadi et al. \cite{AGS10} after the completion of the present study.
They treat the finite element approximation for the Stokes equations with SBCF, using the P1 bubble/P1 element.
Some numerical examples are presented, and an error estimate is announced without a proof.

\section{Settings and results of continuous problems}\label{Sec2}
\subsection{Basic notation}\label{Sec2.1}
Let $\Omega$ be a \emph{polygonal} domain in $\mathbf R^2$.
Throughout this paper, we are concerned with the Stokes equations written in a familiar form
\begin{equation}
	-\nu\Delta u  + \nabla p = f \quad\text{in}\quad \Omega, \qquad \mathrm{div}\, u = 0 \quad\text{in}\quad \Omega, \label{2.9}
\end{equation}
where $\nu>0$ is the viscosity constant; $u$, the velocity field; $p$, the pressure; and $f$, the external force.
As for the boundary, we assume that $\Gamma := \partial\Omega$ is a union of two non-overlapping parts, that is,
\begin{equation*}
	\Gamma = \overline\Gamma_0 \cup \overline\Gamma_1, \quad \Gamma_0 \cap \Gamma_1 = \emptyset,
\end{equation*}
where $\Gamma_0, \Gamma_1$ are relatively nonempty open subsets of $\Gamma$.
Moreover, $\overline\Gamma_1$ is assumed to coincide with \emph{whole one side} of the polygon $\Omega$ for the sake of simplicity.
Two endpoints of the line segment $\overline\Gamma_1$ are respectively denoted by $M_1$ and $M_{m+1}$; the meaning of these subscripts is clarified in Section \ref{Sec3.1}.

We impose the adhesive boundary condition on $\Gamma_0$, namely,
\begin{equation}
	u=0 \quad\text{on}\quad \Gamma_0, \label{2.10}
\end{equation}
whereas on $\Gamma_1$, we impose one (and only one) of the following boundary conditions of friction type:
\begin{equation}
	u_n=0, \qquad |\sigma_\tau|\le g, \qquad \sigma_\tau u_\tau + g|u_\tau|=0, \label{2.1}
\end{equation}
called the \emph{slip boundary condition of friction type} (SBCF), and
\begin{equation}
	u_\tau=0, \qquad |\sigma_n|\le g, \qquad \sigma_n u_n + g|u_n|=0, \label{2.2}
\end{equation}
called the \emph{leak boundary condition of friction type} (LBCF).
The function $g$, called the \emph{modulus of friction}, is assumed to be continuous on $\overline\Gamma_1$ and strictly positive on $\Gamma_1$.

Here, the definitions of the symbols appearing above are as follows:
\begin{gather*}
	n = {}^t(n_1,n_2) = \text{ outer unit normal to the boundary } \Gamma_1, \\
	\tau = {}^t(n_2,-n_1) = \text{ unit tangential vector to the boundary } \Gamma_1, \\
	u_n = u\cdot n = \text{ normal component of $u$ on } \Gamma_1, \\
	u_\tau = u\cdot\tau = \text{ tangential component of $u$ on } \Gamma_1, \\
	e_{ij}(u) = \frac{1}{2}\left( \frac{\partial u_i}{\partial x_j} + \frac{\partial u_j}{\partial x_i} \right) = \text{ component of rate-of-strain tensor } (1\le i,j\le 2), \\
	T_{ij}(u,p) = -p\delta_{ij} \!+\! 2\nu e_{ij}(u) =\! \text{ component of Cauchy stress tensor } (1 \!\le\! i,j \!\le\! 2), \\
	\sigma(u,p) = \Big( \sum_{j=1}^{2}T_{ij}(u,p)n_j \Big)_{i=1,2} = \text{ stress vector defined on }\Gamma_1, \\
	\sigma_n = \sigma_n(u,p) = \text{ normal component of stress vector defined on }\Gamma_1, \\
	\sigma_\tau = \sigma_{\tau}(u) = \text{ tangential component of stress vector defined on }\Gamma_1. \\
\end{gather*}
\begin{Rem}
	(i) $n$ and $\tau$ are constant vectors because $\Gamma_1$ is a segment.
	
	(ii) $\sigma_\tau$ does not depend on $p$, which is verified by a simple calculation.
	
	(iii) In (\ref{2.1}) and (\ref{2.2}), $g$ acts as the threshold of the tangential and normal stress beyond which non-trivial slip and leak on $\Gamma_1$ may occur, respectively.
	This is why the boundary conditions (\ref{2.1}) and (\ref{2.2}) are said to be ``frictional."
\end{Rem}

\subsection{Function spaces}\label{Sec2.2}
\hspace{1mm}We use the usual Lebesgue spaces $L^2(\Omega), L^\infty(\Omega)$ and Sobolev spaces $H^r(\Omega)$ for a nonnegative integer $r$, together with their standard norms and semi-norms
(for a space of vector-valued functions, we write $L^2(\Omega)^2$, and so on). $H^0(\Omega)$ is understood as $L^2(\Omega)$, and $H^1_0(\Omega)$ denotes the closure of $C_0^\infty(\Omega)$ in $H^1(\Omega)$.
We put $Q = L^2(\Omega)$ and
\begin{equation*}
	\stackrel{\;\circ}Q\,\, = L^2_0(\Omega) = \left\{ \psi\in L^2(\Omega) \,\bigg|\, \int_{\Omega}\psi\,dx = 0 \right\}.
\end{equation*}
$H^s(\Omega)$ is also defined for $s\in \mathbf R_+\setminus \mathbf N$ by the norm
\begin{equation*}
	\|\psi\|_{H^s(\Omega)} = \left( \|\psi\|^2_{H^r(\Omega)} + \sum_{|\alpha|=r}\iint_{\Omega\times\Omega}\frac{|\partial^\alpha \psi(x)-\partial^\alpha \psi(y)|^2}{|x-y|^{2+2\theta}}\,dxdy \right),
\end{equation*}
where $\alpha\in\mathbf N^2$ is a multi-index and $s=r+\theta$, where $r\in\mathbf N$, $0<\theta<1$.

We also use the Sobolev space $H^s(\Gamma)$ defined on the boundary $\Gamma$ for $s\ge 0$.
$H^0(\Gamma)$ is understood as $L^2(\Gamma)$, and we put
\begin{equation*}
	L^2_0(\Gamma_1) = \left\{ \eta\in L^2(\Gamma_1) \,\bigg|\, \int_{\Gamma_1}\eta\,ds = 0 \right\},
\end{equation*}
where $ds$ denotes the surface element of $\Gamma$.
The usual trace operator defined on $H^s(\Omega)$ onto $H^{s-1/2}(\Gamma)$ is denoted by $\psi\mapsto \psi|_\Gamma$ for $s>1/2$; however, we simply write $\psi$ instead of $\psi|_\Gamma$ when there is no ambiguity.
Since $n$ and $\tau$ are constant vectors, we immediately obtain the following:
\begin{Lem}\label{Lem2.3}
	Let $s \ge 0$. For every $\phi \in H^s(\Gamma)^2$ satisfying $\phi=0$ on $\Gamma_0$ and $\phi_n=0$ $(\text{resp. }\phi_\tau=0)$ on $\Gamma_1$, it holds that
	\begin{equation*}
		\|\phi\|_{H^s(\Gamma)^2} = \|\phi_\tau\|_{H^s(\Gamma_1)} \qquad (\text{resp. }\|\phi\|_{H^s(\Gamma)^2} = \|\phi_n\|_{H^s(\Gamma_1)}).
	\end{equation*}
\end{Lem}

In addition, we require the so-called Lions-Magene space $H^{\frac12}_{00}(\Gamma_1)$ (see \cite[Section I.11]{LM68}) with its norm defined by
\begin{equation*}
	\|\eta\|_{H^{\frac12}_{00}(\Gamma_1)} = \left( \|\eta\|_{H^{\frac12}(\Gamma_1)}^2 + \int_{\Gamma_1}\frac{|\eta(x)|^2}{\rho(x)}\,ds \right)^\frac{1}{2},
\end{equation*}
where $\rho(x) = \mathrm{dist}(x,\{M_1,M_{m+1}\})$ is the distance from $x\in\Gamma_1$ to the extreme points of $\Gamma_1$ along $\Gamma_1$.
We use this space for only one purpose described in the following lemma.
\begin{Lem}\label{Lem2.1}
The trace operator maps $\{\psi \!\in\! H^1(\Omega) \,|\, \psi \!=\! 0 \text{ on }\Gamma_0\}$ onto $H^{\frac12}_{00}\!(\Gamma_1)$.
\end{Lem}
\begin{proof}
	See \cite[Theorem 1.5.2.3]{G85}.
\end{proof}
\begin{Rem}
	The lemma implies that the extension to $\Gamma$ by zero of an arbitrary function in $H^{\frac12}_{00}(\Gamma_1)$ belongs to $H^{\frac12}(\Gamma)$.
\end{Rem}

Now we let $V=H^1(\Omega)^2$ and introduce the following two closed subspaces of $V$:
\begin{gather}
	V_n = \{ v\in H^1(\Omega)^2 \,|\, v=0\text{ on }\Gamma_0, \quad v_n=0\text{ on }\Gamma_1 \}, \label{2.16} \\
	V_\tau = \{ v\in H^1(\Omega)^2 \,|\, v=0\text{ on }\Gamma_0, \quad v_\tau=0\text{ on }\Gamma_1 \}, \label{2.17}
\end{gather}
which corresponds to the velocity space for SBCF and LBCF, respectively.
Combining the above two lemmas with the usual trace theorem, we see that
\begin{Lem}\label{Lem2.2}
	{\rm(i)} For every $v\in V_n$ $(\text{resp. }v\in V_\tau)$, it holds that
	\begin{align*}
		&v_\tau|_{\Gamma_1}\in H^\frac12_{00}(\Gamma_1) \qquad\text{and}\qquad \|v_\tau\|_{H^\frac12(\Gamma_1)} \le C\|v\|_{H^1(\Omega)^2} \\
		( \text{resp. } &v_n|_{\Gamma_1}\in H^\frac12_{00}(\Gamma_1) \qquad\text{and}\qquad \|v_n\|_{H^\frac12(\Gamma_1)} \le C\|v\|_{H^1(\Omega)^2} ),
	\end{align*}
	with the constant $C$ independent of $v$.
	
	\noindent{\rm(ii)} Every $\eta\in H^\frac12_{00}(\Gamma_1)$ admits an extension $v\in V_n$ $(\text{resp. }v\in V_\tau)$ such that
	\begin{equation*}
		v_\tau=\eta \,\,(\text{resp. }v_n=\eta) \quad\text{on}\quad\Gamma_1 \qquad\text{and}\qquad \|v\|_{H^1(\Omega)^2} \le C\|\eta\|_{H^\frac12(\Gamma_1)},
	\end{equation*}
	with the constant $C$ independent of $\eta$.
\end{Lem}

\subsection{Bilinear forms and barrier terms of friction}\label{Sec2.3}
Let us introduce
\begin{align}
	&a(u,v) = 2\nu\sum_{i,j=1}^{2}\int_{\Omega}e_{ij}(u)e_{ij}(v)\,dx & (u,v\in H^1(\Omega)^2), \\
	&b(v,q) = -\int_{\Omega}\mathrm{div}v\,q\,dx & (v\in H^1(\Omega)^2,\, q\in L^2(\Omega)), \\
	&j(\eta) = \int_{\Gamma_1}g|\eta|\,ds & (\eta\in L^2(\Gamma_1)). \label{2.11}
\end{align}
The bilinear forms $a$ and $b$ are continuous with their operator norms $\|a\|$ and $\|b\|$, respectively, being bounded.
As a readily obtainable consequence of Korn's inequality (\cite[Lemma 6.2]{KO88}), there exists a constant $\alpha > 0$ such that
\begin{equation}
	a(v,v) \ge \alpha\|v\|_{H^1(\Omega)^2}^2 \qquad (\forall v\in V,\, v=0\text{ on }\Gamma_0). \label{2.5}
\end{equation}
This implies that $a$ is coercive on $V_n$ and $V_\tau$.
We simply write $j(v_\tau)$ and $j(v_n)$ to express $j(v_\tau|_{\Gamma_1})$ and $j(v_n|_{\Gamma_1})$, respectively.
Then, $j(v_\tau)$ and $j(v_n)$, called the barrier terms of friction, are continuous functional on $V$ because $j(\eta)$ is bounded on $L^2(\Gamma_1)$.

\subsection{Green's formula}\label{Sec2.4}
For all $(u,p)\in H^{1+\epsilon}(\Omega)^2\times H^\epsilon(\Omega)$ with $\epsilon>1/2$ satisfying ${\rm div}\,u=0$, we obtain Green's formula as follows:
\begin{equation*}
	(-\nu\Delta u + \nabla p, v)_{L^2(\Omega)^2} = a(u,v) + b(v,p) - \int_{\Gamma}\sigma(u,p)\cdot v\,ds \qquad (\forall v\in H^1(\Omega)^2),
\end{equation*}
where the stress vector $\sigma(u,p)$ is defined in Section \ref{Sec2.1}.
In fact, the line integral over $\Gamma$ appearing in the right-hand side is well defined because $\sigma(u,p)\in H^{\epsilon-\frac12}(\Gamma) \subset L^2(\Gamma)$.
However, if we have only a lower regularity, say $(u,p)\in H^1(\Omega)^2\times L^2(\Omega)$, then the definition of $\sigma(u,p)$ in Section \ref{Sec2.1} becomes ambiguous.
We thus propose a redefinition of $\sigma(u,p)$ as a functional on $H^{\frac12}(\Gamma)$ as follows.
\begin{Def}\label{Def2.1}
	Let $(u,p)\in H^1(\Omega)^2\times L^2(\Omega)$ with $\mathrm{div}\,u=0$. When $-\nu\Delta u + \nabla p$ is represented by $f\in L^2(\Omega)^2$ in the distribution sense, that is,
	\begin{equation*}
		\left< -\nu\Delta u + \nabla p, v \right>_{H^1_0(\Omega)^2} = (f,v)_{L^2(\Omega)^2} \qquad (\forall v\in H^1_0(\Omega)^2),
	\end{equation*}
	we define $\sigma(u,p) \in(H^{\frac12}(\Gamma)^2)'$ by
	\begin{equation}
		\left< \sigma(u,p), v \right>_{H^{\frac12}(\Gamma)^2} = a(u,v) + b(v,p) - (f, v)_{L^2(\Omega)^2} \qquad (v\in H^1(\Omega)^2). \label{2.8}
	\end{equation}
	Here and hereafter, for a Banach space $X$, we denote the dual space of $X$ by $X'$ and the duality pairing between $X$ and $X'$ by $\left<\cdot, \cdot\right>_X$.
\end{Def}
\begin{Rem}
	The functional $\sigma(u,p)$ is well defined according to the trace theorem and the fact that the right-hand side of (\ref{2.8}) vanishes if $v=0$ on $\Gamma$, i.e., $v\in H^1_0(\Omega)^2$.
	In addition, this definition of $\sigma(u,p)$ agrees with the previous one if $u$ and $p$ are sufficiently smooth to belong to $H^{1+\epsilon}(\Omega)^2\times H^\epsilon(\Omega)$ with $\epsilon>1/2$.
\end{Rem}
In particular, we see that $\sigma_\tau$ and $\sigma_n$ are characterized in $H^{\frac12}_{00}(\Gamma_1)'$ by
\begin{align*}
	&\langle\sigma_\tau, v_\tau\rangle_{H^{\frac12}_{00}(\Gamma_1)} = a(u,v) + b(v,p) - (f,v)_{L^2(\Omega)^2} &(v\in V_n), \\
	&\langle\sigma_n, v_n\rangle_{H^{\frac12}_{00}(\Gamma_1)} = a(u,v) + b(v,p) - (f,v)_{L^2(\Omega)^2} &(v\in V_\tau),
\end{align*}
respectively, in view of Lemma \ref{Lem2.2}.

\subsection{Variational formulation to the Stokes problem with SBCF}\label{Sec2.5}
Herein we assume $f\in L^2(\Omega)^2$ and $g\in C^1(\overline\Gamma_1)$ with $g>0$ on $\Gamma_1$.
With $V_n$ defined by (\ref{2.16}) and $\stackrel{\;\circ}Q\;=L^2_0(\Omega)$, we introduce a weak formulation of the Stokes equations (\ref{2.9}) under (\ref{2.10}) and (\ref{2.1}) as follows.
\begin{Prob}[PDE]
	Find $(u,p)\in V_n \;\times \stackrel{\;\circ}Q$ such that $\sigma_\tau=\sigma_\tau(u)$ is well defined and the slip boundary condition of friction type $(\ref{2.1})$ is satisfied, that is,
	\begin{numcases}{}
		a(u,v) + b(v,p) - (\sigma_\tau,v_\tau)_{L^2(\Gamma_1)} = (f,v)_{L^2(\Omega)^2} & $(\forall v\in V_n)$, \hspace{1cm} \label{2.4} \\[-1.5mm]
		b(u,q) = 0 & $(\forall q\in\; \stackrel{\;\circ}Q)$, \hspace{1cm} \\
		\sigma_\tau/g\in L^\infty(\Gamma_1) \qquad\text{and}\qquad |\sigma_\tau|\le g \quad\text{a.e. on }\Gamma_1, \label{2.14} \\
		\sigma_\tau u_\tau + g|u_\tau| = 0 \quad\text{a.e. on } \Gamma_1. \label{2.12}
	\end{numcases}
	Note that $\sigma_\tau\in L^2(\Gamma_1)$ follows from $(\ref{2.14})$, and thus $(\ref{2.4})$ makes sense.
\end{Prob}

Another formulation by a variational inequality proposed in \cite{F94} is
\begin{Prob}[VI]
	Find $(u,p)\in V_n \;\times \stackrel{\;\circ}Q$ such that
	\begin{numcases}{\hspace{-1cm}}
		\!\!a(u,v-u) + b(v-u,p) + j(v_\tau) - j(u_\tau) \ge (f,v-u)_{L^2(\Omega)^2} & $(\forall v \!\in\! V_n)$, \label{2.3} \\[-1mm]
		\!\!b(u,q) = 0 & $(\forall q \!\in\, \stackrel{\;\circ}Q)$.
	\end{numcases}
\end{Prob}

The following theorem concerning the existence and uniqueness is essentially derived from \cite[Theorems 2.1--2.3]{F94}.
\begin{Thm}\label{Thm2.1}
	{\rm (i)} Problems {\rm PDE} and {\rm VI} are equivalent in the sense that $(u,p)\in V_n\times\! \stackrel{\;\circ}Q$ solves Problem {\rm PDE} if and only if it solves Problem {\rm VI}.
	
	{\rm (ii)} Problem {\rm VI} has a unique solution.
\end{Thm}
\begin{Rem}
	In \cite{F94}, another definition of $\sigma_\tau = \nu \frac{\partial u_\tau}{\partial n}$ is employed and it is supposed that $\Gamma$ is smooth, with $\overline\Gamma_0 \cap \overline\Gamma_1\neq\emptyset$.
	However, some slight modification, which is not essential, makes the proofs in \cite{F94} applicable to our own situation.
\end{Rem}

\subsection{Variational formulation to the Stokes problem with LBCF}\label{Sec2.6}
As in the previous subsection, using $V_\tau$ defined by (\ref{2.17}) and $Q=L^2(\Omega)$, we introduce a weak formulation of the Stokes equations (\ref{2.9}) under (\ref{2.10}) and (\ref{2.2}) as follows.
\begin{Prob}[PDE]
	Find $(u,p)\in V_\tau\times Q$ such that $\sigma_n=\sigma_n(u,p)$ is well defined and the leak boundary condition of friction type $(\ref{2.2})$ is satisfied, that is,
	\begin{numcases}{}
		a(u,v) + b(v,p) - (\sigma_n,v_n)_{L^2(\Gamma_1)} = (f,v)_{L^2(\Omega)^2} & $(\forall v\in V_\tau)$, \hspace{1cm} \label{2.13} \\
		b(u,q) = 0 & $(\forall q\in Q)$, \hspace{1cm} \\
		\sigma_n/g\in L^\infty(\Gamma_1) \qquad\text{and}\qquad |\sigma_n|\le g \quad\text{a.e. on}\quad \Gamma_1, \label{2.15} \\
		\sigma_n u_n + g|u_n| = 0 \quad\text{a.e. on}\quad \Gamma_1.
	\end{numcases}
	Note that $\sigma_n\in L^2(\Gamma_1)$ follows from $(\ref{2.15})$, and thus $(\ref{2.13})$ makes sense.
\end{Prob}

Another formulation by a variational inequality proposed in \cite{F94} is
\begin{Prob}[VI]
	Find $(u,p)\in V_\tau\times Q$ such that
	\begin{numcases}{\hspace{-1cm}}
		\!\!a(u,v-u) + b(v-u,p) + j(v_n) - j(u_n) \ge (f,v-u)_{L^2(\Omega)^2} & $(\forall v \!\in\! V_\tau)$, \\
		\!\!b(u,q) = 0 & $(\forall q \!\in\! Q)$.
	\end{numcases}
\end{Prob}

We recall the existence and (non)uniqueness theorem derived from \cite[Theorems 3.1--3.3 and Remark 3.2]{F94}.
\begin{Thm}\label{Thm2.2}
	{\rm (i)} Problems {\rm PDE} and {\rm VI} are equivalent in the sense that $(u,p)\in V_\tau\times Q$ solves Problem {\rm PDE} if and only if it solves Problem {\rm VI}.

	{\rm (ii)} Problem {\rm VI} has at least one solution, the velocity part of which is unique.
	
	{\rm (iii)} If $(u,p)$ and $(u,p^*)$ are two solutions of Problem {\rm VI} $($therefore, Problem {\rm PDE}$)$, there exists a unique constant $\delta\in\mathbf R$ such that
	\begin{equation*}
		p = p^* + \delta \qquad\text{and}\qquad \sigma_n(u,p) = \sigma_n(u,p^*) - \delta.
	\end{equation*}
	
	{\rm (iv)} Under the assumptions in {\rm (iii)}, if we suppose $u_n\neq0$ on $\Gamma_1$, then $\delta=0$. Namely, a solution of Problem {\rm VI} is unique.
\end{Thm}
\begin{Rem}
	Although the definition of $\sigma_n = -p + \nu \frac{\partial u_n}{\partial n}$ and the hypotheses on the boundary in \cite{F94} are apparently different from ours,
	we can complete the proof with only a non-essential modification of the original one in \cite{F94}.
\end{Rem}

\section{Finite element approximation}\label{Sec3}
\subsection{Triangulation}\label{Sec3.1}
Let $\{\mathscr{T}_h\}_h$ be a sequence of triangulations of a polygon $\Omega$, where $h$ denotes the length of the greatest side.
As usual, we assume that
\begin{itemize}
	\item $T_i\cap T_j$ is a side, a node, or $\emptyset$ for all $T_i,T_j\in\mathscr{T}_h(i\neq j)$.
	\item $\displaystyle{\bigcup_{i}T_i = \overline{\Omega}}$ and the boundary vertices belong to $\Gamma$.
	\item When $h$ tends to $0$, each triangle in $\mathscr{T}_h$ contains a circle of radius $Kh$ and it is contained in a circle of radius $K'h$ for some constants $K,K'>0$ independent of $h$.
	\item Each triangle has at least one vertex that is not on $\Gamma$.
\end{itemize}
The one-dimensional meshes of $\Gamma$ and $\overline\Gamma_1$ inherited from the triangulation $\mathscr T_h$ are denoted respectively by $\mathscr E_h$ and $\mathscr E_h|_{\overline\Gamma_1}$.
For the sets of nodes, we use
\begin{align*}
	\Sigma_h' &= \text{set of all vertices of triangles in }\mathscr{T}_h, \\
	\Sigma_h'' &= \text{set of all midpoints of sides of triangles in }\mathscr{T}_h, \\
	\Sigma_h &= \Sigma_h'\cup\Sigma_h'', \\
	\Gamma_{0,h} &= \overline\Gamma_0\cap\Sigma_h, \\
	\Gamma_{1,h} &= \overline\Gamma_1\cap\Sigma_h = \{M_1, M_{3/2}, M_2, \cdots, M_m, M_{m+1/2}, M_{m+1}\}, \\
	\stackrel\circ\Gamma_{1,h} &= \Gamma_1\cap\Sigma_h = \Gamma_{1,h}\setminus \{M_1, M_{m+1}\},
\end{align*}
where the subscripts of $M_i$'s are numbered such that
\begin{itemize}
	\item $M_i$'s, for $i=1,2,\cdots,m+1$, are all vertices of triangles in $\mathscr T_h$, which are located in $\overline\Gamma_1$ and are arranged in ascending order along $\overline\Gamma_1$.
	\item $M_{i+1/2}$ is the midpoint of $M_i$ and $M_{i+1}$ for $i=1,2,\cdots,m$.
\end{itemize}
In particular, $\Gamma_{0,h}\cap\Gamma_{1,h} = \overline\Gamma_{0}\cap\overline\Gamma_{1} = \{M_1,M_{m+1}\}$.
We denote each side with endpoints $M_i,M_{i+1}$ by $e_i=[M_i,M_{i+1}]$ and its length by $|e_i|=|M_iM_{i+1}|$, for $i=1,2,\cdots,m$.

\subsection{Approximate function spaces}
We employ the so-called P2/P1 element, defining $V_h\subset V=H^1(\Omega)^2$ and $Q_h\subset Q=L^2(\Omega)$ by
\begin{gather*}
	V_h = \left\{ v_h\in C^0(\overline\Omega)^2 \,\Big|\, v_h|_{T}\in\mathscr{P}_2(T)^2 \quad (\forall T\in\mathscr{T}_h) \right\}, \\
	Q_h = \left\{ q_h\in C^0(\overline\Omega) \,\Big|\, q_h|_{T}\in\mathscr{P}_1(T)^2 \quad (\forall T\in\mathscr{T}_h) \right\},
\end{gather*}
where $\mathscr{P}_k(T)$ denotes the set of all polynomial functions of degree $k$ on $T$ ($k=1,2$).
To approximate $V_n$, $V_\tau$, and $\stackrel{\;\circ}Q\; =L^2_0(\Omega)$, we set
\begin{gather*}
	V_{nh} = \left\{ v_h\in V_h \,\Big|\, v_h(M)=0\;(\forall M\in\Gamma_{0,h}), \quad v_{hn}(M)=0\; (\forall M\in\; \stackrel\circ\Gamma_{1,h}) \right\}, \\
	V_{\tau h} = \left\{ v_h\in V_h \,\Big|\, v_h(M)=0\;(\forall M\in\Gamma_{0,h}), \quad v_{h\tau}(M)=0\; (\forall M\in\; \stackrel\circ\Gamma_{1,h}) \right\}, \\[-1mm]
	\stackrel{\circ~}{Q_h}\; = Q_h \cap L^2_0(\Omega),
\end{gather*}
together with \vspace{-2mm}
\begin{gather*}
	\stackrel{\circ\,}{V_h}\; = V_h\cap H^1_0(\Omega)^2, \\[-1mm]
	V_{nh,\sigma} = \big\{ v_h\in V_h \,\big|\, b(v_h,q_h)=0\, (\forall q_h\in\, \stackrel{\circ~}{Q_h}) \big\}, \\
	V_{\tau h,\sigma} = \big\{ v_h\in V_h \,\big|\, b(v_h,q_h)=0\, (\forall q_h\in Q_h) \big\}.
\end{gather*}
Here, $v_{hn}$ and $v_{h\tau}$ denote $v_h\cdot n$ and $v_h\cdot \tau$, respectively.
By a simple observation we see that $V_{nh}\subset V_n$, $V_{\tau h}\subset V_\tau$, $Q_h\subset Q$, $\stackrel{\circ~}{Q_h} \;\subset\; \stackrel{\;\circ}Q$,  and $\stackrel{\circ\,}{V_h}\; = V_{nh}\cap H^1_0(\Omega)^2 = V_{\tau h}\cap H^1_0(\Omega)^2$.

The quadratic Lagrange interpolation operator $\mathcal I_h: C^0(\overline\Omega)^2\to V_h$ and $L^2$-projection operator $\Pi_h: Q\to Q_h$ are defined in the usual sense, that is,
\begin{gather*}
	\mathcal I_hv\in V_h \quad\text{and}\quad (\mathcal I_hv)(M) = v(M) \qquad (\forall v\in V,\, \forall M\in\Sigma_h), \\
	\Pi_hq\in Q_h \quad\text{and}\quad \int_{\Omega}\left( q - \Pi_hq \right)q_h\,dx = 0 \qquad (\forall q\in Q,\, \forall q_h\in Q_h).
\end{gather*}
It is easy to verify that $\mathcal I_hv\in V_{nh}$ (resp. $\mathcal I_hv\in V_{\tau h}$) if $v\in V_n \cap C^0(\overline\Omega)^2$ (resp. $v\in V_\tau \cap C^0(\overline\Omega)^2$) and that $\Pi_hq\in\; \stackrel{\circ~}{Q_h}$ if $q_h\in\; \stackrel{\;\circ}Q$.
The following results for the interpolation error are standard (for example, see \cite{BS07}) and are used without special emphasis in our error analysis:
\begin{gather}
	\|v - \mathcal I_hv\|_{H^1(\Omega)^2} \le Ch^\epsilon \|v\|_{H^{1+\epsilon}(\Omega)^2} \qquad (\forall v\in H^{1+\epsilon}(\Omega)^2), \label{3.41} \\
	\|q - \Pi_hq\|_{L^2(\Omega)} \le Ch^\epsilon \|v\|_{H^\epsilon(\Omega)} \qquad (\forall q\in H^\epsilon(\Omega)),
\end{gather}
where $0 < \epsilon \le 2$ and the constant $C>0$ depends only on $\Omega$. Note that $H^{1+\epsilon}(\Omega)^2 \subset C^0(\overline\Omega)^2$ for $\epsilon>0$.
Furthermore, the estimate on the boundary, together with Lemma \ref{Lem2.3} and the trace theorem, gives
\begin{gather}
	\|v_\tau - (\mathcal I_hv)_\tau\|_{L^2(\Gamma_1)} \le Ch^{\frac12 + \epsilon}\|v_\tau\|_{H^{\frac12+\epsilon}(\Gamma_1)} \le Ch^{\frac12 + \epsilon}\|v\|_{H^{1+\epsilon}(\Omega)^2} \label{3.42} \\
	\left( \text{resp. } \|v_n - (\mathcal I_hv)_n\|_{L^2(\Gamma_1)} \le Ch^{\frac12 + \epsilon}\|v_n\|_{H^{\frac12+\epsilon}(\Gamma_1)} \le Ch^{\frac12 + \epsilon}\|v\|_{H^{1+\epsilon}(\Omega)^2} \right)
\end{gather}
for all $v\in V_n\cap H^{1+\epsilon}(\Omega)^2$ (resp. $v\in V_\tau\cap H^{1+\epsilon}(\Omega)^2$).

For approximate functions defined on the boundary $\Gamma_1$, we define
\begin{gather*}
	\Lambda_h = \Big\{\mu_h\in C^0(\overline\Gamma_1) \,\Big|\, \mu_h|_e\in\mathscr{P}_2(e)\;(\forall e\in\mathscr E_h|_{\overline\Gamma_1}),\, \mu_h(M_1)\!=\!\mu_h(M_{m+1})\!=\!0 \Big\}, \\
	\tilde{\Lambda}_h = \left\{\mu_h\in\Lambda_h \,\Big|\, |\mu_h(M)|\le 1\;(\forall M\in\; \stackrel\circ\Gamma_{1,h}) \right\}.
\end{gather*}
By a simple calculation, we find that (see also Lemma \ref{Lem3.2}(i))
\begin{equation*}
	\Lambda_h = \left\{ v_{h\tau}|_{\Gamma_1} \,\Big|\, v_h\in V_{nh} \right\} = \left\{ v_{hn}|_{\Gamma_1} \,\Big|\, v_h\in V_{\tau h} \right\}. \\
\end{equation*}
The space $\Lambda_h$ also becomes a Hilbert space if we define its inner product by
\begin{align}
	&(\lambda_h,\mu_h)_{\Lambda_h} \!=\! \frac16\sum_{i=1}^{m}|e_i|\Big( g_i\lambda_{h,i}\mu_{h,i} + 4g_{i+\frac12}\lambda_{h,i+\frac12}\mu_{h,i+\frac12} + g_{i+1}\lambda_{h,i+1}\mu_{h,i+1} \Big) \notag \\[-2mm]
	&\hspace{8cm}	(\lambda_h,\mu_h \in\Lambda_h), \label{3.43}
\end{align}
which approximates ${\int_{\Gamma_1}\!g\mu_h\lambda_h\,ds}$ by Simpson's formula.
Here and in what follows, we occasionally write $g_i,\lambda_{h,i+\frac12},\cdots$ instead of $g(M_i),\lambda_h(M_{i+\frac12}),\cdots$, and so on.
Since $g$ is assumed to be positive on $\Gamma_1$ (particularly, on $\stackrel\circ\Gamma_{1,h}$), $(\cdot,\cdot)_{\Lambda_h}$ is indeed positive definite.
Let us denote the projection operator from the Hilbert space $\Lambda_h$ onto its closed convex subset $\tilde\Lambda_h$ by $\mathrm{Proj}_{\tilde\Lambda_h}$.
It is explicitly expressed as
\begin{equation}\label{3.55}
	\mathrm{Proj}_{\tilde\Lambda_h}(\mu_h)(M) =
	\begin{cases}
		 +1 & \qquad\text{if}\quad \mu_h(M) > 1 \\
		 \mu_h(M) & \qquad\text{if}\quad |\mu_h(M)| \le 1 \\
		 -1 & \qquad\text{if}\quad \mu_h(M) < -1
	\end{cases}
	\qquad (\forall M\in \Gamma_{1,h}),
\end{equation}
for each $\mu_h\in \Lambda_h$.

Finally, to approximate $j$ given in (\ref{2.11}), we introduce $j_h$ as
\begin{align}
	j_h(\eta_h) = \frac16\sum_{i=1}^{m}|e_i|\Big( g_i|\eta_{h,i}| + 4g_{i+\frac12}|\eta_{h,i+\frac12}| + g_{i+1}|\eta_{h,i+1}| \Big) \qquad (\eta_h\in\Lambda_h), \label{3.44}
\end{align}
again with Simpson's formula in mind. Clearly, $j_h$ is a positive, continuous, and positively homogeneous functional defined on $\Lambda_h$.
This definition of $j_h$ is motivated by \cite[Section IV.2.6]{GLT81} and \cite[Section II.5.4]{G84}.

\subsection{Inf-sup conditions}\label{Sec3.3}
Hereafter, we denote various constants independent of $h$ by $C$ and those depending on $h$ by $C(h)$, unless otherwise stated.
In this subsection, two types of inf-sup conditions concerning the approximate spaces of the velocity and pressure are considered.
The first one is the ``$H^1_0$-$L^2_0$" type and well known, while the second one is the ``$H^1$-$L^2$" type and seems to be new.
\begin{Lem}\label{Lem3.1}
	{\rm (i)} There exists a constant $\beta > 0$ independent of $h$ such that
	\begin{equation*}
		\beta\|q_h\|_{L^2(\Omega)} \le \sup_{v_h\in \stackrel{\circ\,}{V_h}}\frac{b(v_h,q_h)}{\|v_h\|_{H^1(\Omega)^2}} \qquad (\forall q_h\in\; \stackrel{\circ~}{Q_h}).
	\end{equation*}
	
	{\rm (ii)} Let $f_1$ and $f_2$ be functions in $L^2(\Omega)^2$ and $L^2(\Omega)$, respectively.
	Then there exists a unique $(u_h,p_h)\in \stackrel{\circ\,}{V_h} \times\stackrel{\circ~}{Q_h}$ such that
	\begin{numcases}{}
		a(u_h,v_h) + b(v_h,p_h) = (f_1,v_h)_{L^2(\Omega)^2} & $(\forall v_h\in\; \stackrel{\circ\,}{V_h})$, \notag \\[-1mm]
		b(u_h,q_h) = (f_2,q_h)_{L^2(\Omega)} & $(\forall q_h\in\; \stackrel{\circ~}{Q_h})$. \notag
	\end{numcases}
	Moreover, $(u_h,p_h)$ satisfies
	\begin{equation*}
		\|u_h\|_{H^1(\Omega)^2} + \|p_h\|_{L^2(\Omega)} \le C(\|f_1\|_{L^2(\Omega)^2} + \|f_2\|_{L^2(\Omega)}),
	\end{equation*}
	where the constant $C$ depends only on $\|a\|, \alpha, \beta$.
\end{Lem}
\begin{proof}
	See \cite[Chapter 12]{BS07}. 
\end{proof}
\begin{Rem}\label{Rem3.1}
	Since $\stackrel{\circ\,}{V_h} \;\subset V_{nh}$, we immediately deduce from (i) that
	\begin{equation*}
		\beta\|q_h\|_{L^2(\Omega)} \le \sup_{v_h\in V_{nh}}\frac{b(v_h,q_h)}{\|v_h\|_{H^1(\Omega)^2}} \qquad (\forall q_h\in\; \stackrel{\circ~}{Q_h}).
	\end{equation*}
\end{Rem}

\begin{Lem}\label{Lem4.2}
	There exists a constant $\beta > 0$ independent of $h$ such that
	\begin{equation*}
		\beta\|q_h\|_{L^2(\Omega)} \le \sup_{v_h\in V_{\tau h}}\frac{b(v_h,q_h)}{\|v_h\|_{H^1(\Omega)^2}} \qquad (\forall q_h\in Q_h).
	\end{equation*}
\end{Lem}
\begin{proof}
	Let us take an arbitrary $p_h \in Q_h$ and define $\eta_h\in\Lambda_h$ by %
	\begin{equation*}
		\eta_h(M_i) =
		\begin{cases}
			\displaystyle{\frac{i-1}{m/2}\xi} \qquad & (i=\frac22,\frac32,\cdots,\frac{m+2}2) \\
			\displaystyle{\frac{m-i+1}{m/2}\xi} \qquad & (i = \frac{m+2}2, \frac{m+3}2, \cdots, \frac{2m+2}2),
		\end{cases}
	\end{equation*}
	where $\xi := -\frac2{|\Gamma_1|}(p_h,1)_{L^2(\Omega)}$ ($|\Gamma_1|$ denotes the length of $\Gamma_1$).
	According to Lemma \ref{Lem3.2}(i), which is preceded by this lemma only for the sake of convenience, we can choose $\hat u_h \in V_{\tau h}$ such that $\hat u_{hn} = \eta_h$ on $\Gamma_1$ and
	\begin{equation}
		\|\hat u_h\|_{H^1(\Omega)^2} \le C\|\eta_h\|_{H^{1/2}(\Gamma_1)}. \label{4.1}
	\end{equation}
	Then, by direct computation we deduce that
	\begin{equation}
		\int_{\Gamma_1}\hat u_{hn}\,ds = \int_{\Gamma_1}\eta_h\,ds = \frac{|\Gamma_1|}2 \xi = -(p_h,1)_{L^2(\Omega)}, \label{4.5}
	\end{equation}
	and that
	\begin{equation*}
		\|\eta_h\|_{H^1(\Gamma_1)}^2 = \left( \frac{|\Gamma_1|}{3} + \frac{4}{|\Gamma_1|} \right)\xi^2 \le C|(p_h,1)_{L^2(\Omega)}|^2 \le C\|p_h\|_{L^2(\Omega)}^2.
	\end{equation*}
	The latter estimate implies
	\begin{equation}
		\|\eta_h\|_{H^\frac12(\Gamma_1)} \le C\|\eta_h\|_{H^1(\Gamma_1)} \le C\|p_h\|_{L^2(\Omega)}. \label{4.2}
	\end{equation}
	From (\ref{4.1}) and (\ref{4.2}), we have
	\begin{equation}
		\|\hat u_h\|_{H^1(\Omega)^2} \le C\|p_h\|_{L^2(\Omega)}. \label{4.34}
	\end{equation}
	
	For $\hat u_h$ constructed above, it follows from Lemma \ref{Lem3.1}(ii) that there exists a unique $(u_h^*,p_h^*) \in\; \stackrel{\circ\,}{V_h} \times \stackrel{\circ~}{Q_h}$ such that
	\begin{numcases}{}
		a(u_h^*,v_h) + b(v_h,p_h^*) = 0 & $(\forall v_h\in\; \stackrel{\circ\,}{V_h})$, \\
		b(u_h^*,q_h) = (p_h,q_h)_{L^2(\Omega)} - b(\hat u_h,q_h) & $(\forall q_h\in\; \stackrel{\circ~}{Q_h})$, \label{4.4}
	\end{numcases}
	together with the estimate
	\begin{align}
		\|u_h^*\|_{H^1(\Omega)^2} &\le C\|p_h + \mathrm{div}\,\hat u_h\|_{L^2(\Omega)} \le C(\|p_h\|_{L^2(\Omega)} + \|\hat u_h\|_{H^1(\Omega)^2}) \notag \\
			&\le C\|p_h\|_{L^2(\Omega)}. \label{4.3}
	\end{align}
	Here we have used (\ref{4.34}) to derive (\ref{4.3}).
	
	Now, setting $u_h = \hat u_h + u_h^* \in V_{\tau h}$ and decomposing $p_h$ as $p_h = p_h^0 + \delta_h$, where $p_h^0 := p_h - (p_h,1)_{L^2(\Omega)^2}/|\Omega| \in\; \stackrel{\circ~}{Q_h}$
	\,and\, $\delta_h := (p_h,1)_{L^2(\Omega)^2}/|\Omega|$ ($|\Omega|$ denotes the area of $\Omega$), we see from (\ref{4.4}) and (\ref{4.5}) that
	\begin{align}
		b(u_h,p_h) &= b(u_h,p_h^0) + b(u_h,\delta_h) = (p_h,p_h^0)_{L^2(\Omega)} - \delta_h\int_\Omega \mathrm{div}u_h\,dx \notag \\
			&= (p_h,p_h^0)_{L^2(\Omega)} - \delta_h\int_{\Gamma_1}\hat u_{hn}\,ds = (p_h,p_h^0)_{L^2(\Omega)} + (p_h,\delta)_{L^2(\Omega)} \notag \\
			&= \|p_h\|_{L^2(\Omega)}^2. \label{4.6}
	\end{align}
	On the other hand, it follows from (\ref{4.34}) and (\ref{4.3}) that
	\begin{equation}
		\|u_h\|_{H^1(\Omega)^2} \le C\|p_h\|_{L^2(\Omega)}. \label{4.7}
	\end{equation}
	From (\ref{4.6}) and (\ref{4.7}), we conclude
	\begin{equation*}
		\sup_{v_h\in V_h}\frac{b(v_h,p_h)}{\|v_h\|_{H^1(\Omega)^2}} \ge \frac{b(u_h,p_h)}{\|u_h\|_{H^1(\Omega)^2}} \ge C\|p_h\|_{L^2(\Omega)}.
	\end{equation*}
	This completes the proof.
\end{proof}
\begin{Rem}
	We can regard this result as a discrete analogue of \cite[Lemma 2.2]{S04}.
\end{Rem}

\subsection{Discrete extension theorems}\label{Sec3.4}
Let us investigate some discrete extensions of functions given on the boundary $\Gamma_1$ to that defined on the whole domain $\Omega$.
\begin{Lem}\label{Lem3.2}
	{\rm (i)} Every $\eta_h\in\Lambda_h$ admits an extension $u_h\in V_{nh}$ $(\text{resp. }u_h\in V_{\tau h})$ such that
	\begin{equation}
		u_{h\tau} = \eta_h ~(\text{resp. }u_{hn} = \eta_h) \text{ on }\Gamma_1 \quad\text{and }\quad \|u_h\|_{H^1(\Omega)^2} \le C\|\eta_h\|_{H^{\frac12}(\Gamma_1)}\, . \label{3.1}
	\end{equation}
	
	{\rm (ii)} For all $\eta_h\in\Lambda_h$ $(\text{resp. }\;\eta_h\in \Lambda_h\cap L^2_0(\Gamma_1))$, we can choose $u_h$ in {\rm (i)} such that
	\begin{equation*}
		u_h\in V_{nh,\sigma} \qquad \big( \text{resp. }u_h\in V_{\tau h,\sigma} \big).
	\end{equation*}
\end{Lem}
\begin{proof}
	(i) Let $\eta_h\in\Lambda_h$.
	We discuss only the construction of $u_h\in V_{nh}$, because we can construct $u_h\in V_{\tau h}$ in a similar manner by replacing $n$ with $\tau$ and vice versa.
	Define a piecewise quadratic polynomial $\phi_h$ on $\Gamma$ by
	\begin{equation*}
		\begin{cases}
			\phi_h|_e \in\mathscr P_2(e)^2 & (\forall e\in \mathscr E_h), \\
			\phi_h(M)=0 & (\forall M\in\Gamma_{0,h}), \\
			\phi_{hn}(M)=0 \text{ and } \phi_{h\tau}(M)=\eta_h(M) & (\forall M\in\; \stackrel\circ\Gamma_{1,h}).
		\end{cases}
	\end{equation*}
	We find that
	\begin{equation}
		\phi_h=0 \text{ on }\Gamma_0, \qquad \phi_{hn}=0 \text{ on }\Gamma_1, \qquad\text{and}\qquad \phi_{h\tau}=\eta_h \text{ on }\Gamma_1, \label{3.2}
	\end{equation}
	and thus we obtain
	\begin{equation}
		\|\phi_h\|_{H^{\frac12}(\Gamma)^2} = \|\eta_h\|_{H^\frac12(\Gamma_1)} \label{3.53}
	\end{equation}
	with the aid of Lemma \ref{Lem2.3}.
	
	Now according to the property of the discrete lifting operator (see \cite[Theorem 5.1]{BG98}), there exists $u_h\in V_h$ satisfying
	\begin{equation}
		u_h = \phi_h\text{ on }\Gamma \qquad\text{and}\qquad \|u_h\|_{H^1(\Omega)^2} \le C\|\phi_h\|_{H^\frac12(\Gamma)^2}. \label{3.3}
	\end{equation}
	We conclude $u_h\in V_{nh}$ and (\ref{3.1}) from (\ref{3.2})--(\ref{3.3}).
	
	(ii) First, take an arbitrary $\eta_h\in\Lambda_h$ and consider an extension to $V_{nh,\sigma}$.
	It follows from (i) that there exists $\hat u_h\in V_{nh}$ such that $\hat u_{h\tau}=\eta_h$ on $\Gamma_1$ and
	\begin{equation}
		\|\hat u_h\|_{H^1(\Omega)^2} \le C\|\eta_h\|_{H^\frac12(\Gamma_1)}. \label{3.4}
	\end{equation}
	For such $\hat u_h$, by Lemma \ref{Lem3.1}, we can find $(u_h^*,p_h^*)\in\; \stackrel{\circ\,}{V_h}\times\stackrel{\circ~}{Q_h}$ satisfying
	\begin{numcases}{}
		a(u_h^*,v_h) + b(v_h,p_h^*) = 0 & $(\forall v_h\in\; \stackrel{\circ\,}{V_h})$, \\
		b(u_h^*,q_h) = -b(\hat u_h,q_h) & $(\forall q_h\in\; \stackrel{\circ~}{Q_h})$, \label{3.5}
	\end{numcases}
	together with the estimate
	\begin{equation}
		\|u_h^*\|_{H^1(\Omega)^2} \le C\|\mathrm{div}\,\hat u_h\|_{L^2(\Omega)} \le C\|\eta_h\|_{H^\frac12(\Gamma_1)}, \label{3.6}
	\end{equation}
	where the last inequality holds from (\ref{3.4}).
	Now, choosing $u_h = u_h^* + \hat u_h\in V_{nh}$, we deduce that $u_h\in V_{nh,\sigma}$ from (\ref{3.5}), that $u_{h\tau}=\hat u_{h\tau}=\eta_h$ because $u_h^*\in\stackrel{\circ\,}{V_h}$,
	and that $\|u_h\|_{H^1(\Omega)^2} \le C\|\eta_h\|_{H^\frac12(\Gamma_1)}$ from (\ref{3.4}) and (\ref{3.6}).
	
	Next, we let $\eta_h\in \Lambda_h \cap L^2_0(\Gamma_1)$ and construct $u_h\in V_{\tau h}$ in the same manner as above by replacing $n$ with $\tau$ and vice versa.
	Then, it remains only to show that $b(u_h,1)=0$ because we already know that $b(u_h,q_h) = 0$ if $q_h\in\; \stackrel{\circ~}{Q_h}$. We can verify $b(u_h,1)=0$ as
	\begin{equation*}
		b(u_h,1) = -\int_\Omega \mathrm{div}u_h\,dx = -\int_{\Gamma_1} u_{hn}\,ds = -\int_{\Gamma_1} \eta_h\,ds = 0.
	\end{equation*}
	This completes the proof. 
\end{proof}

\subsection{Properties of $(\cdot,\cdot)_{\Lambda_h}$ and $j_h$}\label{Sec3.5}
Let us establish several relationships between the inner product of $\Lambda_h$ and the functional $j_h$, given by (\ref{3.43}) and (\ref{3.44}), respectively.
We use a signature function $\mathrm{sgn}(x)$ in the usual sense defined by
\begin{equation*}
	\mathrm{sgn}(x) =
	\begin{cases}
		1 & (x > 0) \\
		0 & (x = 0) \\
		-1 & (x < 0).
	\end{cases}
\end{equation*}
\begin{Lem}\label{Lem3.3}
	{\rm (i)} If $u_h \in V_{nh}$ $(\text{resp. } u_h \in V_{\tau h})$ and $\lambda_h \in \tilde{\Lambda}_h$, then
	\begin{equation*}
		(u_{h\tau}, \lambda_h)_{\Lambda_h} \le j_h(u_{h\tau}) \qquad (\text{resp.}\quad (u_{hn}, \lambda_h)_{\Lambda_h} \le j_h(u_{hn})).
	\end{equation*}
	
	{\rm (ii)} Under the assumptions of {\rm (i)}, the following properties are equivalent:
	\begin{indentation}{30pt}{0pt}\noindent
		{\rm (a)} $(u_{h\tau}, \lambda_h)_{\Lambda_h} \ge j_h(u_{h\tau})$ \qquad $(\text{resp. } (u_{hn}, \lambda_h)_{\Lambda_h} \ge j_h(u_{hn}))$. \\
		{\rm (b)} $(u_{h\tau}, \lambda_h)_{\Lambda_h} = j_h(u_{h\tau})$ \qquad $(\text{resp. } (u_{hn}, \lambda_h)_{\Lambda_h} = j_h(u_{hn}))$. \\
		{\rm (c)} $(u_{h\tau}, \mu_h - \lambda_h)_{\Lambda_h} \le 0$ \qquad $(\text{resp. } (u_{hn}, \mu_h - \lambda_h)_{\Lambda_h} \le 0)$ \qquad $(\forall\mu_h \in \tilde\Lambda_h)$. \\
		{\rm (d)} If $M \in \stackrel\circ\Gamma_{1,h}$ and $u_{h\tau}(M) \neq 0$ $(\text{resp. } u_{hn}(M) \neq 0)$, then
		\begin{equation*}
			\lambda_h(M) = \mathrm{sgn}(u_{h\tau}(M)) \qquad (\text{resp. } \lambda_h(M) = \mathrm{sgn}(u_{hn}(M))).
		\end{equation*}
		{\rm (e)} $\lambda_h = \mathrm{Proj}_{\tilde{\Lambda}_h}(\lambda_h + \rho u_{h\tau})$ \quad $(\text{resp. } \lambda_h = \mathrm{Proj}_{\tilde{\Lambda}_h}(\lambda_h + \rho u_{hn})$ \qquad $(\forall\rho \ge 0)$.
	\end{indentation}
	
	{\rm (iii)} When $\lambda_h\in\Lambda_h$, the following properties are equivalent:
	\begin{indentation}{30pt}{0pt}\noindent
		{\rm (a)} $\lambda_h\in\tilde\Lambda_h$. \\
		{\rm (b)} $(\eta_h,\lambda_h)_{\Lambda_h} \le j_h(\eta_h) \qquad (\forall\eta_h\in\Lambda_h)$.
	\end{indentation}
	
	{\rm (iv)} When $\lambda_h\in\Lambda_h$, the following properties are equivalent:
	\begin{indentation}{30pt}{0pt}\noindent
		{\rm (a)} $(\eta_h, \lambda_h)_{\Lambda_h} = 0 \qquad (\forall\eta_h \in \Lambda_h\cap L^2_0(\Gamma_1))$. \\
		{\rm (b)} There exists a unique constant $\delta_h\in\mathbf R$ such that
			\begin{equation*}
				\lambda_h(M) = \frac{\delta_h}{g(M)} \qquad (\forall M\in \stackrel\circ\Gamma_{1,h}).
			\end{equation*}
	\end{indentation}
\end{Lem}

\begin{proof}
	We establish statements (i) and (ii) only for the case $u_h \in V_{nh}$, because the proof remains valid when $u_h\in V_{\tau h}$, with $n$ replaced by $\tau$ and vice versa.
	
	(i) This is obvious because $|\lambda_h(M)| \le 1$ for all $M\in\Gamma_{1,h}$ if $\lambda_h\in\tilde\Lambda_h$.
	
	(ii) (a)$\Rightarrow$(b) Since we have already proved the converse inequality in (i), statement (b) immediately follows from (a).
	
	(b)$\Rightarrow$(c) Let (b) be valid. From (i), it holds that
	\begin{equation*}
		(u_{h\tau},\mu_h-\lambda_h)_{\Lambda_h} = (u_{h\tau},\mu_h)_{\Lambda_h} - j_h(u_{h\tau}) \le 0 \qquad (\forall\mu_h\in\tilde\Lambda_h).
	\end{equation*}
	
	(c)$\Rightarrow$(d) Assume that (c) is valid and consider an arbitrary $M\in\; \stackrel\circ\Gamma_{1,h}$ such that $u_{h\tau}(M)\neq0$.
	Let us define $\mu_h\in\tilde\Lambda_h$ by
	\begin{equation*}
		\mu_h(N) =
		\begin{cases}
			\lambda_h(N) & \text{if}\quad N\in\; \stackrel\circ\Gamma_{1,h}\setminus\{M\} \\
			\mathrm{sgn}(u_{h\tau}(M)) & \text{if}\quad N=M .
		\end{cases}
	\end{equation*}
	When $M\in\Sigma_h'$, we can write $M=M_i$ for some $1<i<m+1$. Now, by assumption we have
	\begin{equation*}
		(u_{h\tau},\mu_h-\lambda_h)_{\Lambda_h} = \frac{g(M)}6\Big(|e_{i-1}|+|e_i|\Big) \Big(|u_{h\tau}(M)| - \lambda_{h}(M)u_{h\tau}(M)\Big) \le 0.
	\end{equation*}
	This implies that $\lambda_h(M)=\mathrm{sgn}(u_{h\tau}(M))$ because $|\lambda_h(M)| \le 1$ and $u_{h\tau}(M)\neq 0$.
	Similarly, when $M\in\Sigma_h''$, we can write $M=M_{i+\frac12}$ for some $1\le i\le m$.
	Then, by assumption we obtain
	\begin{equation*}
		(u_{h\tau},\mu_h-\lambda_h)_{\Lambda_h} = \frac23 g(M)|e_i|\Big(|u_{h\tau}(M)| - \lambda_{h}(M)u_{h\tau}(M)\Big) \le 0,
	\end{equation*}
	from which $\lambda_h(M)=\mathrm{sgn}(u_{h\tau}(M))$ follows.
	
	(d)$\Rightarrow$(a) If (d) is true, then we see that
	\begin{equation*}
		(u_{h\tau},\lambda_h)_{\Lambda_h} = \frac16\sum_{i=1}^m |e_i|\Big(g_i|u_{h\tau,i}| + 4g_{i+\frac12}|u_{h\tau,i+\frac12}| + g_{i+1}|u_{h\tau,i+1}|\Big) = j_h(u_{h\tau}).
	\end{equation*}
	
	(c)$\Leftrightarrow$(e) This is a direct consequence of a general property of projection operators. In fact, we obtain
	\begin{align*}
		 &(u_{h\tau},\mu_h-\lambda_h)_{\Lambda_h} \le 0 &(\forall \mu_h\in\tilde\Lambda_h) \\
		\iff &(\lambda_h + \rho u_{h\tau} - \lambda_h, \mu_h - \lambda_h)_{\Lambda_h} \le 0 &(\forall \mu_h\in\tilde\Lambda_h, \forall\rho\ge 0) \\
		\iff &\lambda_h = \mathrm{Proj}_{\tilde\Lambda_h}(\lambda_h + \rho u_{h\tau}) &(\forall\rho\ge 0).
	\end{align*}
	
	(iii) (a)$\Rightarrow$(b) This is already shown in (i).
	
	(b)$\Rightarrow$(a) Let (b) be valid and consider an arbitrary $M\in\; \stackrel\circ\Gamma_{1,h}$.
	Define $\eta_h\in \Lambda_h$ by
	\begin{equation*}
		\eta_h(N) =
		\begin{cases}
			0 &\text{if}\quad N\in\; \stackrel\circ\Gamma_{1,h}\setminus\{M\} \\
			+1 \text{ or }-1 &\text{if}\quad N=M .
		\end{cases}
	\end{equation*}
	When $M\in\Sigma_h'$, we can write $M=M_i$ for some $1<i<m+1$. By assumption, we obtain $(\eta_h,\lambda_h)_{\Lambda_h} \le j_h(\eta_h)$, which leads to
	\begin{equation*}
		\frac16 \Big(|e_{i-1}| + |e_i|\Big)g(M)\Big(\lambda_h(M) \pm 1\Big) \le 0.
	\end{equation*}
	This implies that $|\lambda_h(M)| \le 1$.
	We obtain the same result when $M\in\Sigma_h''$ in a similar way.
	Therefore, we conclude that $\lambda_h\in\tilde\Lambda_h$.
	
	(iv) (b)$\Rightarrow$(a) Let (b) be valid and consider such $\delta_h$.
	Because Simpson's formula is exact for quadratic polynomials,	for all $\eta_h \in \Lambda_h\cap L^2_0(\Gamma_1)$, we have
	\begin{align*}
		(\eta_h, \lambda_h)_{\Lambda_h} &= \sum_{i=1}^m \frac{|e_i|}6 \left( g_i\eta_{h,i}\lambda_{h,i} + 4g_{i+\frac12}\eta_{h,i+\frac12}\lambda_{h,i+\frac12} + g_{i+1}\eta_{h,i+1}\lambda_{h,i+1} \right) \\
			&= \delta_h \sum_{i=1}^m\frac{|e_i|}6 \left( \eta_{h,i} + 4\eta_{h,i+\frac12} + \eta_{h,i+1} \right) \\
			&= \delta_h \int_{\Gamma_1}\eta_h\,ds = 0.
	\end{align*}
	
	(a)$\Rightarrow$(b) Let (a) be valid and consider $i\in \{2,3,\cdots,m\}$. Let us make $\eta_h\in\Lambda_h$ vanish except on $e_{i-1}\cup e_i$.
	Then, statement (a) is equivalently written as

	(a$'$)\quad{\it
		For all $\eta_{h,i-\frac12},\eta_{h,i},\eta_{h,i+\frac12} \in \mathbf R$, if we assume
		\begin{equation*}
			|e_{i-1}|(4\eta_{h,i-\frac12} + \eta_{h,i}) + |e_i|(\eta_{h,i} + 4\eta_{h,i+\frac12}) = 0,
		\end{equation*}
		then we have
		\begin{equation*}
			|e_{i-1}|\Big\{4\eta_{h,i-\frac12}(g\lambda_h)_{i-\frac12} + \eta_{h,i}(g\lambda_h)_i\Big\} + |e_i|\Big\{\eta_{h,i}(g\lambda_h)_i + 4\eta_{h,i+\frac12}(g\lambda_h)_{i+\frac12}\Big\} = 0.
		\end{equation*}
	}
	
	Now, if we take $\eta_h|_{e_{i-1}\cup e_i}$ such that $\eta_{h,i}=0$ and $4\eta_{h,i-\frac12}|e_{i-1}| = -4\eta_{h,i+\frac12}|e_i| = 1$,
	it follows from (a$'$) that $(g\lambda_h)_{i-\frac12} = (g\lambda_h)_{i+\frac12}$.
	Similarly, if we take $\eta_h|_{e_{i-1}\cup e_i}$ such that $\eta_{h,i-\frac12}=0$ and $\eta_{h,i}(|e_{i-1}| + |e_i|) = -4\eta_{h,i+\frac12}|e_i| = 1$,
	it follows again from (a$'$) that $(g\lambda_h)_i = (g\lambda_h)_{i+\frac12}$. Hence, $(g\lambda_h)_{i-\frac12} = (g\lambda_h)_i = (g\lambda_h)_{i+\frac12}$.
	
	Repeating the above procedure for $i=2,3,\cdots,m$, we conclude that there exists $\delta_h\in\mathbf R$ such that
	\begin{equation*}
		g(M)\lambda_h(M) = \delta_h \qquad (\forall M\!\in\, \stackrel\circ\Gamma_{1,h}).
	\end{equation*}
	This completes the proof.
\end{proof}

The following mesh-dependent inf-sup condition is important to deduce the unique existence of the Lagrange multiplier $\lambda_h\in\Lambda_h$, which appears in Sections \ref{Sec4} and \ref{Sec5}.
\begin{Lem}
	There exists a positive constant $\beta_h$ depending on $h$ such that
	\begin{gather}
		\beta_h\|\eta_h\|_{\Lambda_h} \le \sup_{v_h \in V_{nh}} \frac{(v_{h\tau},\eta_h)_{\Lambda_h}}{\|v_h\|_{H^1(\Omega)^2}} \qquad (\forall\eta_h \in \Lambda_h) \\
		\left( \text{ resp. } \beta_h\|\eta_h\|_{\Lambda_h} \le \sup_{v_h \in V_{\tau h}} \frac{(v_{hn},\eta_h)_{\Lambda_h}}{\|v_h\|_{H^1(\Omega)^2}} \qquad (\forall\eta_h \in \Lambda_h) \right). \label{3.40}
	\end{gather}
\end{Lem}
\begin{proof}
	Because both $\|\cdot\|_{H^\frac12(\Gamma_1)}$ and $\|\cdot\|_{\Lambda_h}$ are norms defined on $\Lambda_h$, which is of a finite dimension, they are equivalent.
	Hence there exists a constant $C(h)$ depending on $h$ such that
	\begin{equation*}
		\|\eta_h\|_{H^\frac12(\Gamma_1)} \le C(h)\|\eta_h\|_{\Lambda_h} \qquad (\forall \eta_h\in\Lambda_h).
	\end{equation*}
	
	Now, we let $\eta_h\in\Lambda_h$ and choose $u_h\in V_{nh}$ satisfying (\ref{3.1}). Then we conclude
	\begin{align*}
		\sup_{v_h\in V_{nh}}\frac{(v_{h\tau},\eta_h)_{\Lambda_h}}{\|v_h\|_{H^1(\Omega)^2}} &\ge \frac{(u_{h\tau},\eta_h)_{\Lambda_h}}{\|u_h\|_{H^1(\Omega)^2}}
			= \frac{\|\eta_h\|^2_{\Lambda_h}}{\|u_h\|_{H^1(\Omega)^2}} \ge C(h)\frac{\|\eta_h\|_{H^\frac12(\Gamma_1)}}{\|u_h\|_{H^1(\Omega)^2}}\|\eta_h\|_{\Lambda_h} \notag \\
			&\ge C(h)\|\eta_h\|_{\Lambda_h}.
	\end{align*}
	We can obtain (\ref{3.40}) in a similar way. 
\end{proof}

\subsection{Error between $j$ and $j_h$}\label{Sec3.6}
We begin with some generalization of \cite[Lemma IV.1.3]{GLT81} concerning the error between $j$ and $j_h$, which is necessary later in our error analysis.
\begin{Lem}\label{Lem3.4}
	{\rm (i)} There hold
	\begin{gather}
		j_h(\eta_h) \le C\|\eta_h\|_{L^2(\Gamma_1)} \qquad (\forall \eta_h\in\Lambda_h), \label{3.45} \\
		\|\eta_h\|_{\Lambda_h} \le C\|\eta_h\|_{L^2(\Gamma_1)} \qquad (\forall \eta_h\in\Lambda_h), \label{3.35}
	\end{gather}
	with the constant $C$ depending only on $g$ and $\Gamma_1$.
	
	{\rm (ii)} If $g\in C^1(\overline\Gamma_1)$, then for all $0\le s\le 1$, we have
	\begin{equation}
		|j_h(\eta_h) - j(\eta_h)| \le Ch^s\|\eta_h\|_{H^s(\Gamma_1)} \qquad (\forall \eta_h\in\Lambda_h), \label{3.21}
	\end{equation}
	with the constant $C$ depending only on $g$ and $\Gamma_1$.
\end{Lem}
\begin{proof}
	(i) Let $\eta_h\in\Lambda_h$.
	On each segment $e_i = [M_i,M_{i+1}]$, take two points denoted by $M_{i+\frac16}$ and $M_{i+\frac56}$, whose meaning is understood naturally, for $i=1,2,\cdots,m$.
	Let us define a piecewise constant function $r_h(g\eta_h)$ on $\overline\Gamma_1$ by
	\begin{align}
		 &r_h(g\eta_h) \notag \\
		= &\sum_{i=1}^m\Big\{ g_i\eta_{h,i}\chi_{[M_i,M_{i+\frac16}]} + g_{i+\frac12}\eta_{h,i+\frac12}\chi_{[M_{i+\frac16},M_{i+\frac56}]} + g_{i+1}\eta_{h,i+1}\chi_{[M_{i+\frac56},M_{i+1}]} \Big\}, \label{hhhhhh}
	\end{align}
	where $\chi_A$ denotes the characteristic function of $A\subset\overline\Gamma_1$.

	Then we have
	\begin{equation}
		j_h(\eta_h) = \int_{\Gamma_1}\!|r_h(g\eta_h)|\,ds \le |\Gamma_1|^{1/2}\|r_h(g\eta_h)\|_{L^2(\Gamma_1)}. \label{3.51}
	\end{equation}
	By direct computation, it follows that
	\begin{align}
		 \|r_h(g\eta_h)\|_{L^2(\Gamma_1)}^2 &= \sum_{i=1}^m\frac{|e_i|}6\left( g_i^2\eta_{h,i}^2 + 4g_{i+\frac12}^2\eta_{h,i+\frac12}^2 + g_{i+1}^2\eta_{h,i+1}^2 \right) \notag \\
			&\le \sup|g|^2 \sum_{i=1}^m\frac{|e_i|}6\left( \eta_{h,i}^2 + 4\eta_{h,i+\frac12}^2 + \eta_{h,i+1}^2 \right) \notag \\
			&\le \frac52 \sup|g|^2 \sum_{i=1}^m\frac{|e_i|}{15}\Big\{ 2(\eta_{h,i}^2+\eta_{h,i+1}^2) + 8\eta_{h,i+\frac12}^2 \notag \\
			&\hspace{3cm}	- \eta_{h,i}\,\eta_{h,i+1} + 2\eta_{h,i+\frac12}(\eta_{h,i}+\eta_{h,i+1}) \Big\} \notag \\
			&= \frac52 \sup|g|^2\, \|\eta_h\|_{L^2(\Gamma_1)}^2. \label{3.52}
	\end{align}
	Here we have used the inequality
	\begin{align}
		&x^2+4y^2+z^2 \le 2(x^2+z^2)+8y^2-xz+2y(x+z) \label{3.46} \\
		\iff &4\left(y+\frac{x+z}4\right)^2 + \frac34(x-z)^2 \ge 0 \qquad(\forall x,y,z\in \mathbf R) \notag
	\end{align}
	to derive the third line. We conclude (\ref{3.45}) from (\ref{3.51}) and (\ref{3.52}).

	The estimate (\ref{3.35}) follows similarly if we remark that
	\begin{equation*}
		\|\eta_h\|_{\Lambda_h}^2 = \sum_{i=1}^m\frac{|e_i|}6 \left( g_i\eta_{h,i}^2 + 4g_{i+\frac12}\eta_{h,i+\frac12}^2 + g_{i+1}\eta_{h,i+1}^2 \right).
	\end{equation*}
	
	(ii) Let $\eta_h\in\Lambda_h$. First, from the proof of (i), we see that
	\begin{align}
		|j_h(\eta_h) - j(\eta_h)| &\le \int_{\Gamma_1}\!|r_h(g\eta_h) - g\eta_h|\,ds \notag \\
			&\le |\Gamma_1|^{1/2}\|r_h(g\eta_h) - g\eta_h\|_{L^2(\Gamma_1)} \label{cccccc} \\
			&\le |\Gamma_1|^{1/2} \left(\|r_h(g\eta_h)\|_{L^2(\Gamma_1)} + \|g\eta_h\|_{L^2(\Gamma_1)}\right) \notag \\
			&\le C\|\eta_h\|_{L^2(\Gamma_1)}. \label{ffffff}
	\end{align}
	
	Before giving an estimate of $\|r_h(g\eta_h) - g\eta_h\|_{L^2(\Gamma_1)}$ which involves $\|\eta_h\|_{H^1(\Gamma_1)}$, it should be noted that if $\phi \in \mathscr{P}_2(\mathbf{R})$ we have
	\begin{equation*}
		\phi(x) - \phi(a) = (x - a)\phi'\left(\frac{a+x}{2}\right),
	\end{equation*}
	so that
	\begin{align}
		\int_a^b|\phi(x) - \phi(a)|^2\,dx &= \int_a^b|x - a|^2 \left|\phi'\left(\frac{a+x}{2}\right)\right|\,dx \notag \\
			&= 8\int_a^{(a+b)/2}|t - a|^2|\phi'(t)|^2\,dt \notag \\
			&\le 2(b - a)^2\int_a^b|\phi'(t)|^2\,dt. \label{bbbbbb}
	\end{align}
	In view of the Taylor expansion of $g$, we apply (\ref{bbbbbb}) to deduce
	\begin{align}
		I_1(e_i) &:= \int_{M_i}^{M_{i+\frac16}}\!|g_i\eta_{h,i} - g\eta_h|^2\,ds = \int_{M_i}^{M_{i+\frac16}}\!|g_i\eta_{h,i}-g_i\eta_h+g_i\eta_h-g\eta_h|^2\,ds \notag \\
			&\le 2\int_{M_i}^{M_{i+\frac16}}\!g_i^2(\eta_{h,i} - \eta_h)^2\,ds + 2\int_{M_i}^{M_{i+\frac16}}\!(g_i - g)^2\eta_h^2\,ds \notag \\
			&\le 2\sup|g|^2\cdot 2\left( \frac{|e_i|}6 \right)^2\int_{M_i}^{M_{i+\frac16}}\!|\eta_h'|^2\,ds + 2\left( \frac{|e_i|}6\sup|g'| \right)^2\int_{M_i}^{M_{i+\frac16}}\!|\eta_h|^2\,ds \notag \\
			&\le \max\left\{ \frac{\sup|g|^2}9, \frac{\sup|g'|^2}{18} \right\} |e_i|^2 \|\eta_h\|_{H^1([M_i,M_{i+\frac16}])}^2 \notag \\
			&\le Ch^2\|\eta_h\|_{H^1(e_i)}^2, \label{dddddd}
	\end{align}
	for $i=1,2,\cdots,m$. 
	By a similar discussion, we have
	\begin{align}
		I_2(e_i) &:= \int_{M_{i+\frac16}}^{M_{i+\frac12}}\!|g_{i+\frac12}\eta_{h,i+\frac12} - g\eta_h|^2\,ds \le Ch^2\|\eta_h\|_{H^1(e_i)}^2, \\
		I_3(e_i) &:= \int_{M_{i+\frac12}}^{M_{i+\frac56}}\!|g_{i+\frac12}\eta_{h,i+\frac12} - g\eta_h|^2\,ds \le Ch^2\|\eta_h\|_{H^1(e_i)}^2, \\
		I_4(e_i) &:= \int_{M_{i+\frac56}}^{M_{i+1}}\!|g_{i+1}\eta_{h,i+1} - g\eta_h|^2\,ds \le Ch^2\|\eta_h\|_{H^1(e_i)}^2, \label{eeeeee}
	\end{align}
	for each $i$. Therefore, it follows from (\ref{cccccc}) and (\ref{dddddd})--(\ref{eeeeee}) that
	\begin{align*}
		|j_h(\eta_h) - j(\eta_h)|^2 &\le C \|r_h(g\eta_h) - g\eta_h\|_{L^2(\Gamma_1)}^2 \\
			&= C\sum_{i=1}^m \Big( I_1(e_i) + I_2(e_i) + I_3(e_i) + I_4(e_i) \Big) \\
			&\le Ch^2\sum_{i=1}^m \|\eta_h\|_{H^1(e_i)}^2 = Ch^2\|\eta_h\|_{H^1(\Gamma_1)}^2,
	\end{align*}
	so that
	\begin{equation}
		|j_h(\eta_h) - j(\eta_h)| \le Ch\|\eta_h\|_{H^1(\Gamma_1)}. \label{gggggg}
	\end{equation}
	As a consequence of (\ref{ffffff}) and (\ref{gggggg}), we obtain the desired inequality (\ref{3.21}) by Hilbertian interpolation (see \cite[Chapter 14]{BS07}) between $L^2(\Gamma_1)$ and $H^1(\Gamma_1)$. 
\end{proof}

As will be shown in Theorems \ref{Thm3.2} and \ref{Thm4.3} below, the leading term of the error is that between $j_h$ and $j$, which is estimated by (\ref{3.21}) with $s=1/2$.
However, under some additional conditions, we can obtain a sharper estimate than (\ref{3.21}).
\begin{Def}\label{Def3.2}
	An element $\eta_h\in\Lambda_h$ is said to \emph{have a constant sign on every side} if, for any $i=1,2,\cdots,m$, either of the following conditions is satisfied:
	
	\centering {\rm(a)} $\eta_h|_{e_i} \ge 0$\qquad or\qquad {\rm(b)} $\eta_h|_{e_i} \le 0$.
\end{Def}
\begin{Rem}
	Let $\eta_h\in\Lambda_h$ have a constant sign on every side. If $\eta_h\ge0$ on $e_{i-1}$ and $\eta_h\le0$ on $e_i$ for some $i$, then $\eta_h(M_i) = 0$.
\end{Rem}

\begin{Lem}\label{Lem3.5}
	Let $g\in C^2(\overline\Gamma_1)$. If $\eta_h\in\Lambda_h$ has a constant sign on every side, then
	\begin{equation*}
		|j_h(\eta_h) - j(\eta_h)| \le Ch^2\|\eta_h\|_{L^2(\Gamma_1)}.
	\end{equation*}
	Moreover, if $g$ is a polynomial of degree $\le 1$, then $j_h(\eta_h)$ is exact, that is,
	\begin{equation}
		j_h(\eta_h) = j(\eta_h). \label{oooooo}
	\end{equation}
\end{Lem}
\begin{proof}
	Let $\eta_h\in\Lambda_h$ have a constant sign on every side. Because $\eta_h\ge0$ or $\eta_h\le0$ on $e_i$ for each $i=1,2,\cdots,m$ and $g$ is positive on $\Gamma_1$, we have
	\begin{gather*}
		\int_{e_i}\left| r_h(g\eta_h) \right|\,ds = \left| \int_{e_i}r_h(g\eta_h)\,ds \right|, \\
		\int_{e_i}g|\eta_h|\,ds = \left| \int_{e_i}g\eta_h\,ds \right|,
	\end{gather*}
	where $r_h(g\eta_h)$ is defined as (\ref{hhhhhh}). Summing up these terms, we obtain
	\begin{gather*}
		j_h(\eta_h) = \int_{\Gamma_1}\!|r_h(g\eta_h)|\,ds = \sum_{i=1}^m\int_{e_i}\left| r_h(g\eta_h) \right|\,ds = \sum_{i=1}^m\left|\int_{e_i}\!r_h(g\eta_h)\,ds\right|, \\
		j(\eta_h) = \int_{\Gamma_1}\!g|\eta_h|\,ds = \sum_{i=1}^m\int_{e_i}g|\eta_h|\,ds = \sum_{i=1}^m\left|\int_{e_i}\!g\eta_h\,ds\right|.
	\end{gather*}
	Consequently, it follows that
	\begin{equation}
		|j_h(\eta_h) - j(\eta_h)| \le \sum_{i=1}^m\left| \int_{e_i}\!\big(r_h(g\eta_h) - g\eta_h\big)ds \right|. \label{iiiiii}
	\end{equation}
	Let $g_h$ denote the linear Lagrange interpolation of $g$ using the nodes in $\Sigma_h'\cap\Gamma_{1,h}$.
	Namely, $g_h$ is continuous on $\overline\Gamma_1$ and affine on each side $e_i=[M_i,M_{i+1}]$, satisfying $g_h(M_i) = g(M_i)$ for $i=1,2,\cdots,m$.
	Then the Taylor expansion of $g$ implies
	\begin{equation}
		|g_h(x) - g(x)| \le \frac{h^2}8 \sup|g''| \qquad (\forall x\in \overline\Gamma_1). \label{kkkkkk}
	\end{equation}
	
	Now, let us estimate each term appearing in the summation on the right-hand side of (\ref{iiiiii}) by
	\begin{equation}
		\left|\int_{e_i}\!\big(r_h(g\eta_h) - g_h\eta_h\big)ds\right| + \int_{e_i}\!|g_h-g||\eta_h|\,ds. \label{jjjjjj}
	\end{equation}
	Since Simpson's formula is exact for cubic polynomials, we can express
	\begin{align*}
		\int_{e_i}\!g_h\eta_h\,ds &= \frac{|e_i|}6\sum_{i=1}^m\left(g_{h,i}\eta_{h,i} + 4g_{h,i+\frac12}\eta_{h,i+\frac12} + g_{h,i+1}\eta_{h,i+1}\right) \\
			&= \int_{e_i}\!r_h(g\eta_h)\,ds + \frac{2|e_i|}3(g_{h,i+\frac12} - g_{i+\frac12})\eta_{h,i+\frac12}.
	\end{align*}
	Thus, due to (\ref{kkkkkk}), the first term of (\ref{jjjjjj}) is bounded from above by
	\begin{equation*}
		\frac1{12} |e_i|h^2\sup|g''|\, |\eta_{h,i+\frac12}|.
	\end{equation*}
	Note that there holds (cf. (\ref{3.46}))
	\begin{align*}
		&|\eta_{h,i+\frac12}|^2|e_i| \\
		\le &\frac{15}4\cdot \frac{|e_i|}{15}\Big\{ 2(\eta_{h,i}^2+\eta_{h,i+1}^2) + 8\eta_{h,i+\frac12}^2 - \eta_{h,i}\,\eta_{h,i+1} + 2\eta_{h,i+\frac12}(\eta_{h,i}+\eta_{h,i+1}) \Big\} \\
		= &\frac{15}4\|\eta_h\|_{L^2(e_i)}^2 \le 4\|\eta_h\|_{L^2(e_i)}^2
	\end{align*}
	for $i=1,2,\cdots,m$. Then, the sum of the first term of (\ref{jjjjjj}) is estimated as
	\begin{align}
		\sum_{i=1}^m\left|\int_{e_i}\!\big(r_h(g\eta_h) - g_h\eta_h\big)ds\right| &\le \frac1{12} h^2\sup|g''|\sum_{i=1}^m|\eta_{h,i+\frac12}|\,|e_i| \notag \\
			&\le \frac1{12} h^2\sup|g''| \left(\sum_{i=1}^m|\eta_{h,i+\frac12}|^2|e_i|\right)^\frac12 \left(\sum_{i=1}^m|e_i|\right)^\frac12 \notag \\
			&\le \frac1{12} h^2\sup|g''|\cdot 2\|\eta_h\|_{L^2(\Gamma_1)}|\Gamma_1|^{1/2} \notag \\
			&= Ch^2\|\eta_h\|_{L^2(\Gamma_1)}. \label{llllll}
	\end{align}
	Next, the second term of (\ref{jjjjjj}) is estimated by $\frac18h^2\sup|g''|\int_{e_i}\!|\eta_h|\,ds$, which gives
	\begin{align}
		\sum_{i=1}^m\int_{e_i}\!|g_h-g|\,|\eta_h|\,ds &\le \frac18h^2 \sup|g''|\int_{\Gamma_1}\!|\eta_h|\,ds \le \frac18h^2 \sup|g''|\,\|\eta_h\|_{L^2(\Gamma_1)}|\Gamma_1|^{1/2} \notag \\
			&= Ch^2\|\eta_h\|_{L^2(\Gamma_1)}. \label{mmmmmm}
	\end{align}
	Hence we conclude from (\ref{iiiiii}), (\ref{llllll}), and (\ref{mmmmmm}) that
	\begin{equation*}
		|j_h(\eta_h) - j(\eta_h)| \le Ch^2\|\eta_h\|_{L^2(\Gamma_1)}.
	\end{equation*}
	
	If $g$ is a polynomial of degree $\le 1$, then both terms of (\ref{jjjjjj}) vanish because $g_h=g$, from which (\ref{oooooo}) follows.
	This completes the proof.
\end{proof}

\section{Discretization of the Stokes problem with SBCF}\label{Sec4}
\subsection{Existence and uniqueness results}\label{Sec4.1}
We propose approximate problems for Problem VI (therefore, Problem PDE) in the case of SBCF as follows.
\begin{Prob}[VI$_h$]
	Find $(u_h, p_h) \in V_{nh} \;\times \stackrel{\circ~}{Q_h}$ such that
	\begin{numcases}{\hspace{-7mm}}
		\hspace{-2mm} a(u_h,v_h\!-\!u_h) \!+\! b(v_h\!-\!u_h,p_h) \!+\! j_h(\!v_{h\tau}\!) \!-\! j_h(\!u_{h\tau}\!) \!\ge\! (f,v_h\!-\!u_h)\!_{L^2(\Omega)^2} &\hspace{-7.5mm} $(\forall v_h \!\!\in\!\! V_{nh})$\!, \label{3.15} \\[-1mm]
		\hspace{-2mm} b(u_h,q_h) = 0 &\hspace{-7.5mm} $(\forall q_h \!\!\in \stackrel{\circ~}{Q_h})$.
	\end{numcases}
\end{Prob}
\begin{Prob}[VI$_{h,\sigma}$]
	Find $u_h\in V_{nh,\sigma}$ such that
	\begin{equation}
		a(u_h,v_h-u_h) + j_h(v_h) - j_h(u_h) \ge (f,v_h-u_h)_{L^2(\Omega)^2} \qquad (\forall v_h\in V_{nh,\sigma}). \label{3.7}
	\end{equation}
\end{Prob}
\begin{Prob}[VE$_h$]
	Find $(u_h, p_h, \lambda_h) \in V_{nh} \;\times \stackrel{\circ~}{Q_h}\times\tilde\Lambda_h$ such that
	\begin{numcases}{\hspace{-1cm}}
		a(u_h,v_h) + b(v_h,p_h) + (v_{h\tau}, \lambda_h)_{\Lambda_h} = (f,v_h)_{L^2(\Omega)^2} & $(\forall v_h\in V_{nh})$, \label{3.13} \\[-2mm]
		b(u_h,q_h) = 0 & $(\forall q_h\in\; \stackrel{\circ~}{Q_h})$, \label{3.11} \\
		(u_{h\tau},\mu_h - \lambda_h)_{\Lambda_h} \le 0 & $(\forall\mu_h \in \tilde\Lambda_h)$. \label{3.12}
	\end{numcases}
\end{Prob}

Recall that we are assuming $f\in L^2(\Omega)^2$ and $g\in C^0(\overline\Gamma_1)$.
We first establish the existence and uniqueness of these approximate problems.

\begin{Thm}\label{Thm3.1}
	{\rm(i)} Problem $\mathrm{VI}_{h,\sigma}$ admits a unique solution $u_h\in V_{nh,\sigma}$.
	Furthermore, it satisfies the following equation:
	\begin{equation}
		a(u_h,u_h) + j_h(u_{h\tau}) = (f,u_h)_{L^2(\Omega)^2}. \label{3.8}
	\end{equation}
	
	\noindent{\rm(ii)} Problems $\mathrm{VI}_{h,\sigma}$, $\mathrm{VI}_h$, and $\mathrm{VE}_h$ are equivalent in the following sense.

	{\rm(a)} If $u_h \in V_{nh,\sigma}$ is a solution of Problem $\mathrm{VI}_{h,\sigma}$, then there exists a unique $p_h \in\; \stackrel{\circ~}{Q_h}$ such that $(u_h,p_h)$ solves Problem $\mathrm{VI}_h$.
	
	{\rm(b)} If $(u_h,p_h)\in V_{nh} \;\times \stackrel{\circ~}{Q_h}$ is a solution of Problem {\rm VI}$_h$, then there exists a unique $\lambda_h\in \tilde\Lambda_h$
		such that $(u_h,p_h,\lambda_h)$ solves Problem {\rm VE}$_{h}$.
		
	{\rm(c)} If $(u_h,p_h,\lambda_h) \in V_{nh} \;\times\stackrel{\circ~}{Q_h}\times\; \tilde\Lambda_h$ is a solution of Problem {\rm VE}$_h$, then $u_h$ solves Problem {\rm VI}$_{h,\sigma}$.
\end{Thm}
\begin{proof}
	(i) Since the bilinear form $a$ is coercive on $V_{nh}$ and the functional $j_h:V_{nh}\to\mathbf R$ is convex, proper, and lower semi-continuous (actually, continuous) with respect to
	the weak topology, we can apply to Problem VI$_{nh,\sigma}$ a classical existence and uniqueness theorem for second-order elliptic variational inequalities (see \cite[Theorem I.4.1]{G84}).
	Thus, there exists a unique $u_h\in V_{nh,\sigma}$ such that (\ref{3.7}) holds. The equation (\ref{3.8}) follows from (\ref{3.7}) with $v_h=0$ and $2u_h$.

	(ii) (a) Let $u_h \in V_{nh,\sigma}$ be a solution of Problem {\rm VI}$_{h,\sigma}$.
	Taking $u_h \pm v_h$ as a test function in (\ref{3.7}), with an arbitrary $v_h \in\; \stackrel{\circ\,}{V_h}\cap\; V_{nh,\sigma}$, we obtain
	\begin{equation*}
		a(u_h,v_h) = (f,v_h)_{L^2(\Omega)^2} \qquad (\forall v_h \in\; \stackrel{\circ\,}{V_h}\cap\; V_{nh,\sigma}).
	\end{equation*}
	Moreover, from Lemma \ref{Lem3.1}(i), we deduce the unique existence of $p_h \!\in\, \stackrel{\circ~}{Q_h}$ such that
	\begin{equation}
		a(u_h,v_h) + b(v_h,p_h) = (f,v_h)_{L^2(\Omega)^2} \qquad (\forall v_h \in\; \stackrel{\circ\,}{V_h}) \label{3.9}
	\end{equation}
	by a standard argument.
	
	Now we let $v_h \in V_{nh}$ be arbitrary.
	It follows from Lemma \ref{Lem3.2} (ii) that there exists some $w_h \in V_{nh,\sigma}$ such that $w_h = v_h$ on $\Gamma$, which implies
	\begin{equation}
		v_h-w_h \in\; \stackrel{\circ\,}{V_h} \qquad\text{and}\qquad j_h(v_{h\tau})=j_h(w_{h\tau}). \label{3.10}
	\end{equation}
	Since $u_h,w_h\in V_{nh,\sigma}$, we conclude from (\ref{3.7}), (\ref{3.9}), and (\ref{3.10}) that
	\begin{align*}
		& a(u_h,v_h-u_h) + b(v_h-u_h,p_h) + j_h(v_{h\tau}) - j_h(u_{h\tau}) - (f,v_h-u_h)_{L^2(\Omega)^2} \\
		= \,& a(u_h,v_h-w_h) + b(v_h-w_h,p_h) - (f,v_h-w_h)_{L^2(\Omega)^2} \notag \\
			&\qquad\qquad + a(u_h,w_h-u_h) + j_h(w_{h\tau}) - j_h(u_{h\tau}) - (f,w_h-u_h)_{L^2(\Omega)^2} \\
		\ge \,& 0.
	\end{align*}
	Hence $(u_h,p_h)$ is a solution of VI$_h$.
	
	(b) Let $(u_h,p_h)\in V_{nh}\; \times \stackrel{\circ~}{Q_h}$ be a solution of VI$_h$.
	Taking $u_h\pm v_h$ as a test function in (\ref{3.15}), with an arbitrary $v_h\in\; \stackrel{\circ\,}{V_h}$, we have
	\begin{equation}
		a(u_h,v_h) + b(v_h,p_h) = (f,v_h)_{L^2(\Omega)^2} \qquad (\forall v_h \in\; \stackrel{\circ\,}{V_h}). \label{3.19}
	\end{equation}
	Therefore, since $\{v_h\in V_{nh} \,|\, (v_{h\tau}, \eta_h)_{\Lambda_h}=0\; (\forall\eta_h\in\Lambda_h)\} =\; \stackrel{\circ\,}{V_h}$, the inf-sup condition
	given in Lemma \ref{Lem4.2} asserts the unique existence of $\lambda_h \in \Lambda_h$ such that
	\begin{equation}
		a(u_h,v_h) + b(v_h,p_h) + (v_{h\tau}, \lambda_h)_{\Lambda_h} = (f,v_h)_{L^2(\Omega)^2} \qquad (\forall v_h\in V_{nh}). \label{3.16}
	\end{equation}
	Combining (\ref{3.16}) with (\ref{3.15}), we obtain
	\begin{equation}
		(v_{h\tau} - u_{h\tau}, \lambda_h)_{\Lambda_h} \le j_h(v_{h\tau}) - j_h(u_{h\tau}) \qquad (\forall v_h\in V_{nh}), \label{3.17}
	\end{equation}
	which gives, by a triangle inequality, that
	\begin{equation}
		(v_{h\tau} - u_{h\tau}, \lambda_h)_{\Lambda_h} \le j_h(v_{h\tau} - u_{h\tau}) \qquad (\forall v_h\in V_{nh}). \label{3.18}
	\end{equation}
	From (\ref{3.18}) together with Lemma \ref{Lem3.2}(i), we deduce
	\begin{equation}
		(\eta_h, \lambda_h)_{\Lambda_h} \le j_h(\eta_h) \qquad (\forall \eta_h\in \Lambda_h).
	\end{equation}
	Hence Lemma \ref{Lem3.3}(iii) implies that $\lambda_h\in \tilde\Lambda_h$, and (\ref{3.13}) is established.
	It remains only to prove (\ref{3.12}).
	Taking $v_h=0$ in (\ref{3.17}), we have $j_h(u_{h\tau}) \le (u_{h\tau},\lambda_h)_{\Lambda_h}$. 
	This implies (\ref{3.12}) by Lemma \ref{Lem3.3}(ii).
	Therefore, $(u_h,p_h,\lambda_h)$ is a solution of Problem VE$_h$.
	
	(c) Let $(u_h,p_h,\lambda_h) \in V_{nh} \;\times\stackrel{\circ~}{Q_h}\times\; \tilde\Lambda_h$ be a solution of Problem VE$_h$.
	Then we see that $u_h\in V_{nh,\sigma}$ from (\ref{3.11}), and that
	\begin{equation}
		(u_{h\tau}, \lambda_h)_{\Lambda_h} = j_h(u_{h\tau}) \label{3.14}
	\end{equation}
	from (\ref{3.12}) combined with Lemma \ref{Lem3.3}(ii).
	It follows from (\ref{3.13}) and (\ref{3.14}) that
	\begin{align*}
		&a(u_h,v_h-u_h) + j_h(v_{h\tau}) - j_h(u_{h\tau}) - (f,v_h-u_h)_{L^2(\Omega)^2} \notag \\
		= &- b(v_h-u_h, p_h) - (v_{h\tau}-u_{h\tau}, \lambda_h)_{\Lambda_h} + j_h(v_{h\tau}) - j_h(u_{h\tau}) \\
		= &- (v_{h\tau}-u_{h\tau}, \lambda_h)_{\Lambda_h} + j_h(v_{h\tau}) - j_h(u_{h\tau}) \\
		= &j_h(v_h) - (v_{h\tau}, \lambda_h)_{\Lambda_h} \\
		\ge &0 \qquad (\forall v_h\in V_{nh,\sigma}),
	\end{align*}
	where the last inequality holds by Lemma \ref{Lem3.3}(i).
	Hence $u_h$ is a solution of Problem VI$_{h,\sigma}$.
	This completes the proof of Theorem \ref{Thm3.1}. 
\end{proof}

\subsection{Error analysis}\label{Sec4.2}
Before presenting the rate-of-convergence results, we state the following:
\begin{Prop}\label{Prop3.1}
	Let $(u,p)$ be the solution of Problem {\rm VI} and $(u_h,p_h)$, that of Problem {\rm VI}$_h$ for $h>0$. Then,
	
	{\rm(i)} it holds that
	\begin{equation}
		\|u_h\|_{H^1(\Omega)^2} \le \|f\|_{L^2(\Omega)^2}/\alpha. \label{3.20}
	\end{equation}
	
	{\rm(ii)} for every $v_h\in V_{nh}$ and $q_h\in\; \stackrel{\circ~}{Q_h}$, it holds that
	\begin{align}
		&\alpha\|u-u_h\|_{H^1(\Omega)^2}^2 \le a(u-u_h,u-v_h) + b(u_h-u,p-q_h) + b(v_h-u,p_h-p) \notag \\
		&\hspace{2cm} + (\sigma_\tau,v_{h\tau}-u_\tau)_{L^2(\Gamma_1)} + j(u_{h\tau}) - j_h(u_{h\tau}) + j_h(v_{h\tau}) - j(u_\tau). \label{3.22}
	\end{align}
	
	{\rm(iii)} for every $q_h \in\; \stackrel{\circ~}{Q_h}$, it holds that
	\begin{align}
		\|p-p_h\|_{L^2(\Omega)^2} \le \left( 1 + \frac{\|b\|}\beta \right) \|p-q_h\|_{L^2(\Omega)^2} + \frac{\|a\|}\beta \|u-u_h\|_{H^1(\Omega)^2}. \label{3.26}
	\end{align}
\end{Prop}
\begin{proof}
	(i) Since $u_h$ is the solution of Problem VI$_{h,\sigma}$ by Theorem \ref{Thm3.1}(ii), it satisfies (\ref{3.8}).
	Hence Korn's inequality (\ref{2.5}), together with the positiveness of $j_h$, gives
	\begin{align*}
		\alpha\|u_h\|^2_{H^1(\Omega)^2} &\le a(u_h,u_h) = (f,u_h)_{L^2(\Omega)^2} - j_h(u_{h\tau}) \\
			&\le (f,u_h)_{L^2(\Omega)^2} \\
			&\le \|f\|_{L^2(\Omega)^2}\|u_h\|_{L^2(\Omega)^2},
	\end{align*}
	which implies (\ref{3.20}).
	
	(ii) Let $v_h \in V_{nh}$ and $q_h \in Q_h$ be arbitrary. We begin with the following equality:
	\begin{equation*}
		a(u-u_h, u-u_h) = a(u-u_h, u-v_h) - a(u, u_h-u) - a(u_h, v_h-u_h) + a(u, v_h-u).
	\end{equation*}
	We bound from above the second term of the right-hand side by (\ref{2.3}) with $v=u_h$, the third one by (\ref{3.15}) with $v_h$ itself,
	and rewrite the fourth one by (\ref{2.4}) with $v=v_h-u$. Consequently,
	\begin{align*}
		 &a(u-u_h,u-u_h) \notag \\
		\le &a(u-u_h,u-v_h) + b(u_h-u,p) + j(u_{h\tau}) - j(u_\tau) - (f,u_h-u)_{L^2(\Omega)^2} \notag \\
			&\hspace{26.2mm} + b(v_h-u_h,p_h) + j_h(v_{h\tau}) - j_h(u_{h\tau}) - (f,v_h-u_h)_{L^2(\Omega)^2} \notag \\
			&\hspace{26.2mm} - b(v_h-u,p) + (\sigma_\tau,v_{h\tau}-u_\tau)_{L^2(\Gamma_1)} + (f,v_h-u)_{L^2(\Omega)^2} \notag \\
		= &a(u-u_h,u-v_h) + b(u_h-u,p-q_h) + b(v_h-u,p_h-p) \notag \\
			&\hspace{1cm} + (\sigma_\tau,v_{h\tau}-u_\tau)_{L^2(\Gamma_1)} + j(u_{h\tau}) - j_h(u_{h\tau}) + j_h(v_{h\tau}) - j(u_\tau). 
	\end{align*}
	Combining this with Korn's inequality (\ref{2.5}), we conclude (\ref{3.22}).
	
	(iii) Taking $u\pm v$ as a test function in (\ref{2.3}), with an arbitrary $v\in H^1_0(\Omega)^2$, gives
	\begin{equation*}
		a(u,v) + b(v,p) = (f,v)_{L^2(\Omega)^2} \qquad (\forall v\in H^1_0(\Omega)^2).
	\end{equation*}
	On the other hand we know that (\ref{3.19}) holds, and therefore, by subtraction we obtain
	\begin{equation}
		a(u-u_h,v_h) + b(v_h,p-p_h) = 0 \qquad (\forall v_h\in\; \stackrel{\circ\,}{V_h}). \label{3.23}
	\end{equation}
	Now let $q_h \in\, \stackrel{\circ~}{Q_h}$. It is clear that
	\begin{equation}
		\|p-p_h\|_{L^2(\Omega)} \le \|p-q_h\|_{L^2(\Omega)} + \|q_h-p_h\|_{L^2(\Omega)}. \label{3.24}
	\end{equation}
	By Lemma \ref{Lem3.1}(i) together with (\ref{3.23}), we have
	\begin{align}
		\beta\|q_h-p_h\| &\le \sup_{v_h\in \stackrel{\circ\,}{V_h}} \frac{b(v_h,q_h-p_h)}{\|v_h\|_{H^1(\Omega)^2}} = \sup_{v_h\in \stackrel{\circ\,}{V_h}} \frac{b(v_h,q_h-p) + b(v_h,p-p_h)}{\|v_h\|_{H^1(\Omega)^2}} \notag \\
			&= \sup_{v_h\in \stackrel{\circ\,}{V_h}} \frac{b(v_h,q_h-p) - a(u-u_h,v_h)}{\|v_h\|_{H^1(\Omega)^2}} \notag \\
			&\le \|b\|\,\|p-q_h\|_{L^2(\Omega)} + \|a\|\,\|u-u_h\|_{H^1(\Omega)^2}. \label{3.25}
	\end{align}
	The desired inequality (\ref{3.26}) follows from (\ref{3.24}) and (\ref{3.25}). 
\end{proof}

We are now in a position to state the primary result of our error estimates, assuming only the regularity of the exact solution.
\begin{Thm}\label{Thm3.2}
	Let $(u,p)$ be the solution of Problem {\rm VI} and $(u_h,p_h)$ be that of Problem {\rm VI}$_h$\ for $0<h<1$.
	Suppose $g\in C^1(\overline\Gamma_1)$ and $(u,p) \in H^{1+\epsilon}(\Omega)^2 \times H^\epsilon(\Omega)$ with $0< \epsilon\le 2$. Then we have
	\begin{equation}
		\|u - u_h\|_{H^1(\Omega)^2} + \|p - p_h\|_{L^2(\Omega)} \le Ch^{\min\{\epsilon,\frac14\}}. \label{3.30}
	\end{equation}
\end{Thm}
\begin{proof}
	We recall the interpolation error estimates (\ref{3.41})--(\ref{3.42}).
	
	Taking $v_h = \mathcal I_hu, q_h=\Pi_hp$ in (\ref{3.22}) and (\ref{3.26}), we find that
	\begin{align}
		\alpha\|u - u_h\|_{H^1(\Omega)^2}^2 &\le a(u-u_h,u-\mathcal I_hu) + b(u_h-u,p-\Pi_hp) + b(\mathcal I_hu-u,p_h-p) \notag \\
		&\hspace{1cm} + (\sigma_\tau,u_\tau-(\mathcal I_hu)_\tau)_{L^2(\Gamma_1)} + |j(u_{h\tau}) - j_h(u_{h\tau})| \notag \\
		&\hspace{1cm} + |j_h((\mathcal I_hu)_\tau) - j((\mathcal I_hu)_\tau)| + |j((\mathcal I_hu)_\tau) - j(u_\tau)|, \label{3.27}
	\end{align}
	and that
	\begin{align}
		\|p - p_h\|_{L^2(\Omega)} &\le C\left( \|p - \Pi_hp\|_{L^2(\Omega)} + \|u - u_h\|_{H^1(\Omega)^2} \right) \notag \\
			&\le C(h^\epsilon + \|u - u_h\|_{H^1(\Omega)^2}). \label{3.28}
	\end{align}
	Each term of the right-hand side in (\ref{3.27}) is estimated as follows:
	
	1.
	\begin{equation*}
		|a(u-u_h,u-\mathcal I_hu)| \le \|a\|\, \|u-u_h\|_{H^1(\Omega)^2}\|u-\mathcal I_hu\|_{H^1(\Omega)^2} \le Ch^\epsilon \|u-u_h\|_{H^1(\Omega)^2}.
	\end{equation*}
	
	2.
	\begin{equation*}
		|b(u_h-u,p-\Pi_hp)| \le \|b\|\,\|u-u_h\|_{H^1(\Omega)^2}\|p - \Pi_hp\|_{L^2(\Omega)} \le Ch^\epsilon \|u-u_h\|_{H^1(\Omega)^2}.
	\end{equation*}
	
	3. From (\ref{3.28}),
	\begin{align*}
		|b(\mathcal I_hu-u,p_h-p)| &\le \|b\|\,\|\mathcal I_hu-u\|_{H^1(\Omega)^2}\|p_h - p\|_{L^2(\Omega)} \\
			&\le C(h^{2\epsilon} + h^\epsilon\|u - u_h\|_{H^1(\Omega)^2}).
	\end{align*}
	
	4.
	\begin{equation*}
		\left|\left( \sigma_\tau,u_\tau-(\mathcal I_hu)_\tau \right)_{L^2(\Gamma_1)}\right| \le \|\sigma_\tau\|_{L^2(\Gamma_1)} \|u_\tau-(\mathcal I_hu)_\tau\|_{L^2(\Gamma_1)}
			\le Ch^{\frac12 + \epsilon}.
	\end{equation*}
	
	5. By Lemma \ref{Lem3.4}(ii) together with Proposition \ref{Prop3.1}(i),
	\begin{equation*}
		|j(u_{h\tau}) - j_h(u_{h\tau})| \le Ch^\frac12 \|u_{h\tau}\|_{H^\frac12(\Gamma_1)} \le Ch^\frac12 \|u_h\|_{H^1(\Omega)^2} \le Ch^\frac12.
	\end{equation*}
	
	6. Since $\|\mathcal I_hu\|_{H^1(\Omega)^2} \!\le\! \|\mathcal I_hu - u\|_{H^1(\Omega)^2} + \|u\|_{H^1(\Omega)^2} \!\le\! C$, Lemma \ref{Lem3.4}(ii) implies
	\begin{equation*}
		|j((\mathcal I_hu)_\tau) - j_h((\mathcal I_hu)_\tau)| \le Ch^\frac12\|(\mathcal I_hu)_\tau\|_{H^\frac12(\Gamma_1)} \le Ch^\frac12\|\mathcal I_hu\|_{H^1(\Omega)^2} \le Ch^\frac12.
	\end{equation*}
	
	7.
	\begin{equation*}
		|j((\mathcal I_hu)_\tau) - j(u_\tau)| \le \|g\|_{L^2(\Gamma_1)}\|u_\tau-(\mathcal I_hu)_\tau\|_{L^2(\Gamma_1)} \le Ch^{\frac12 + \epsilon}.
	\end{equation*}
	Combining these seven estimates with (\ref{3.27}), we deduce that
	\begin{align*}
		\|u - u_h\|^2_{H^1(\Omega)^2} &\le C\left( h^\epsilon\|u - u_h\|_{H^1(\Omega)^2} + h^{2\epsilon} + h^{\frac12 + \epsilon} + h^\frac12 \right) \notag \\
			&\le C\left( h^{\min\{\epsilon, \frac14\}} \|u - u_h\|_{H^1(\Omega)^2} + h^{\min\{2\epsilon, \frac12\}} \right).
	\end{align*}
	Therefore,
	\begin{equation}
		\|u - u_h\|_{H^1(\Omega)^2} \le Ch^{\min\{\epsilon, \frac14\}}. \label{3.29}
	\end{equation}
	We conclude (\ref{3.30}) from (\ref{3.28}) and (\ref{3.29}) and this completes the proof.
\end{proof}

The previous theorem reveals that the rate of convergence is $O(h^{1/4})$ at best even when the solution is sufficiently smooth.
However, it can be improved if additional conditions about the signs of $u_{h\tau}$ and $(\mathcal I_hu)_\tau$ on $\Gamma_1$ are available.
To formulate the result, we make the following assumptions (Recall Definition \ref{Def3.2} and see Remark \ref{Rem4.1}):

(S1)\quad $(\mathcal I_hu)_\tau$ has a constant sign on every side.

(S2)\quad $u_{h\tau}$ has a constant sign on every side.

(S3)\quad $\mathrm{sgn}(u_\tau) = \mathrm{sgn}((\mathcal I_hu)_\tau)$ on $\Gamma_1$.
\begin{Thm}\label{Thm3.3}
	In addition to the hypotheses in Theorem $\ref{Thm3.2}$, we assume $g\in C^2(\overline\Gamma_1)$ and that $(\mathrm{S}1)$--$(\mathrm{S}3)$ are satisfied. Then we have
	\begin{equation}
		\|u - u_h\|_{H^1(\Omega)^2} + \|p - p_h\|_{L^2(\Omega)} \le Ch^{\min\{\epsilon,1\}}. \label{3.48}
	\end{equation}
	Moreover, if $g$ is a polynomial function of degree $\le 1$, we have
	\begin{equation}
		\|u - u_h\|_{H^1(\Omega)^2} + \|p - p_h\|_{L^2(\Omega)} \le Ch^\epsilon. \label{3.49}
	\end{equation}
\end{Thm}
\begin{proof}
	We first verify that (S3) implies
	\begin{equation}
		\sigma_\tau(\mathcal I_hu)_\tau + g|(\mathcal I_hu)_\tau| = 0 \quad\text{a.e. on }\, \Gamma_1. \label{3.31}
	\end{equation}
	In fact, for each side $e\in\mathscr E_h|_{\Gamma_1}$, if $u_\tau$ vanishes on a subset of $e$ containing more than three points, then the quadratic polynomial
	$(\mathcal I_hu)_\tau$ vanishes on the whole $e$. Otherwise, we have $|u_\tau| > 0$ a.e. on $e$; hence we deduce from (\ref{2.12}), namely,
	\begin{equation}
		\sigma_\tau u_\tau + g|u_\tau| = 0 \quad\text{a.e. on }\, \Gamma_1, \label{3.32}
	\end{equation}
	that $\sigma_\tau = -g\mathrm{sgn}(u_\tau) = -g\mathrm{sgn}((\mathcal I_hu)_\tau)$\, a.e.\, on\, $e$.
	In both cases, it follows that $\sigma_\tau(\mathcal I_hu)_\tau + g|(\mathcal I_hu)_\tau| = 0$\, a.e.\, on\, $e$. Thus (\ref{3.31}) is valid.
	
	It follows from (\ref{3.31}) and (\ref{3.32}) that
	\begin{equation*}
		(\sigma_\tau, (\mathcal I_hu)_\tau - u_\tau)_{L^2(\Gamma_1)} + j((\mathcal I_hu)_\tau) - j(u_\tau) = 0.
	\end{equation*}
	Therefore, taking $v_h = \mathcal I_hu$ and $q_h=\Pi_hp$ in (\ref{3.22}) gives
	\begin{align}
		\alpha\|u - u_h\|_{H^1(\Omega)^2}^2 &\le a(u-u_h,u-\mathcal I_hu) + b(u_h-u,p-\Pi_hp) + b(\mathcal I_hu-u,p_h-p) \notag \\
		&\hspace{1cm} + |j(u_{h\tau}) - j_h(u_{h\tau})| + |j_h((\mathcal I_hu)_\tau) - j((\mathcal I_hu)_\tau)|. \label{3.33}
	\end{align}
	Let us give estimates for each term on the right-hand side.
	We can evaluate the first three terms by the same way as in the proof of Theorem \ref{Thm3.2}.
	By Lemma \ref{Lem3.5}, the fourth and fifth terms are estimated as
	\begin{gather*}
		|j(u_{h\tau}) - j_h(u_{h\tau})| \le Ch^2\|u_{h\tau}\|_{L^2(\Gamma_1)} \le Ch^2\|u_h\|_{H^1(\Omega)^2} \le Ch^2, \\
		|j_h((\mathcal I_hu)_\tau) - j((\mathcal I_hu)_\tau)| \!\le\! Ch^2\|(\mathcal I_hu)_\tau\|_{L^2(\Gamma_1)} \le Ch^2\|\mathcal I_hu\|_{H^1(\Omega)^2} \le Ch^2.
	\end{gather*}
	Consequently, we obtain
	\begin{equation}
		\|u - u_h\|^2_{H^1(\Omega)^2} \le C\left( h^\epsilon\|u - u_h\|_{H^1(\Omega)^2} + h^{2\epsilon} + h^2 \right), \label{3.50}
	\end{equation}
	which leads to
	\begin{equation*}
		\|u - u_h\|_{H^1(\Omega)^2} \le Ch^{\min\{\epsilon, 1\}}.
	\end{equation*}
	The estimate for $\|p - p_h\|_{L^2(\Omega)}$ is similar to the proof of Theorem \ref{Thm3.2}, and then, (\ref{3.48}) follows.
	
	Finally, if $g$ is affine then the fourth and fifth terms in (\ref{3.33}) vanish exactly, according to Lemma \ref{Lem3.5}. Hence we have
	\begin{equation*}
		\|u - u_h\|^2_{H^1(\Omega)^2} \le C\left( h^\epsilon\|u - u_h\|_{H^1(\Omega)^2} + h^{2\epsilon} \right)
	\end{equation*}
	instead of (\ref{3.50}), from which (\ref{3.49}) follows.
\end{proof}
\begin{Rem}\label{Rem4.1}
	Conditions (S1)--(S3) are not so artificial. Assume that $u$, the velocity part of the solution, is continuous on $\overline\Omega$ and that the isolated zeros of
	$u_\tau$ on $\Gamma_1$ are contained in $\Gamma_{1,h}$. If we make $h$ sufficiently small, then we see that (S1) and (S3) are satisfied.
	Therefore, since Theorem \ref{Thm3.2} implies $u_{h\tau}\to u_\tau$ in $H^\frac12(\Gamma_1)$, we can expect (S2) to also be valid; however, its rigorous proof is not easy.
\end{Rem}

\subsection{Numerical realization}\label{Sec4.3}
We propose the following Uzawa-type method to compute the solution of Problem VE$_h$ (therefore, Problem VI$_h$) numerically.

\begin{Alg}\label{Alg3.1}
	Choose an arbitrary $\lambda_h^{(1)} \in \tilde\Lambda_h$ and $\rho > 0$. Iterate the following two steps for $k=1,2,\cdots:$
	
	{\rm\bf \underline{Step 1}}\; With $\lambda_h^{(k)}$ known, determine $(u_h^{(k)},p_h^{(k)}) \in V_{nh} \;\times \stackrel{\circ~}{Q_h}$ by
	\begin{numcases}{\hspace{-5mm}}
		\!\!a(u_h^{(k)},v_h) + b(v_h,p_h^{(k)}) = (f,v_h)_{L^2(\Omega)^2} - (v_{h\tau},\lambda_h^{(k)})_{\Lambda_h} &\hspace{-2mm}$(\forall v_h \!\in\! V_{nh})$, \label{3.34} \\
		\!\!b(u_h^{(k)},q_h) = 0 &$\hspace{-2mm}(\forall q_h\!\in\, \stackrel{\circ~}{Q_h})$. \label{3.54}
	\end{numcases}
	
	{\rm\bf \underline{Step 2}}\; Renew $\lambda_h^{(k+1)} \in \tilde\Lambda_h$ by
	\begin{equation}
		\lambda_h^{(k+1)} = \mathrm{Proj}_{\tilde\Lambda_h}(\lambda_h^{(k)} + \rho u_{h\tau}^{(k)}). \label{3.58}
	\end{equation}
\end{Alg}
\begin{Rem}
	(i) The unique existence of $(u_h^{(k)},p_h^{(k)})$ satisfying (\ref{3.34}) and (\ref{3.54}) is guaranteed by the inf-sup condition mentioned in Remark \ref{Rem3.1}.
	
	(ii) We can regard (\ref{3.58}) as an approximation of
	\begin{equation}
		\lambda_h = \mathrm{Proj}_{\tilde\Lambda_h}(\lambda_h + \rho u_{h\tau}), \label{3.59}
	\end{equation}
	which is equivalent to (\ref{3.12}) by Lemma \ref{Lem3.3}(ii).
\end{Rem}

\begin{Thm}\label{Thm3.4}
	Let $(u_h,p_h,\lambda_h)$ be the solution of Problem $\mathrm{VE}_h$. Under the same notation as Algorithm $\ref{Alg3.1}$,
	there exists a constant $\rho_0>0$ independent of $h$ such that if $\rho$ satisfies $0 < \rho < \rho_0$, then the iterative solution $(u_h^{(k)},p_h^{(k)},\lambda_h^{(k)})$ converges
	to $(u_h,p_h,\lambda_h)$ in $H^1(\Omega)^2\times L^2(\Omega)\times \Lambda_h$, as $k \to \infty$.
\end{Thm}
\begin{proof}
	Subtracting (\ref{3.34}) from (\ref{3.13}) with test functions in $V_{nh,\sigma}$, we obtain
	\begin{equation}
		a(u_h-u_h^{(k)},v_h) + (v_{h\tau}, \lambda_h-\lambda_h^{(k)})_{\Lambda_h} = 0 \qquad (\forall v_h \in V_{nh,\sigma}). \label{3.38}
	\end{equation}
	In particular, we take $v_h=u_h^{(k)}-u_h \in V_{nh,\sigma}$ and apply Korn's inequality (\ref{2.5}) to obtain
	\begin{equation}
		(u_{h\tau}^{(k)}-u_{h\tau}, \lambda_h^{(k)}-\lambda_h)_{\Lambda_h} = -a(u_h^{(k)}-u_h,u_h^{(k)}-u_h) \le -\alpha\|u_h^{(k)}-u_h\|_{H^1(\Omega)^2}^2. \label{3.60}
	\end{equation}	
	Next, we note that $\mathrm{Proj}_{\tilde\Lambda_h}$ given in (\ref{3.55}) satisfies
	\begin{equation}
		\|\mathrm{Proj}_{\tilde\Lambda_h}(\mu_h) - \mathrm{Proj}_{\tilde\Lambda_h}(\eta_h)\|_{\Lambda_h} \le \|\mu_h - \eta_h\|_{\Lambda_h} \qquad (\forall \mu_h,\eta_h \in \Lambda_h), \label{3.56}
	\end{equation}
	as a result of a general property of a projection operator.
	It follows from (\ref{3.56}) with $\mu_h=\lambda_h^{(k)} + \rho u_{h\tau}^{(k)}$ and $\eta_h=\lambda_h + \rho u_{h\tau}$, (\ref{3.58}), (\ref{3.59}), and (\ref{3.60}) that
	\begin{align*}
		\|\lambda_h^{(k+1)}-\lambda_h\|_{\Lambda_h}^2 &\le \|\lambda_h^{(k)}-\lambda_h + \rho (u_{h\tau}^{(k)}-u_{h\tau}) \|_{\Lambda_h}^2 \notag \\
			&= \|\lambda_h^{(k)}-\lambda_h\|_{\Lambda_h}^2 \!+ 2\rho (u_{h\tau}^{(k)}-u_{h\tau}, \lambda_h^{(k)}-\lambda_h)_{\Lambda_h} + \rho^2\|u_{h\tau}^{(k)}-u_{h\tau}\|_{\Lambda_h}^2 \notag \\
			&\le \|\lambda_h^{(k)}-\lambda_h\|_{\Lambda_h}^2 - 2\alpha\rho\|u_h^{(k)}-u_h\|_{H^1(\Omega)^2}^2 + \rho^2\|u_{h\tau}^{(k)}-u_{h\tau}\|_{\Lambda_h}^2.
	\end{align*}
	Therefore, since $\|u_{h\tau}^{(k)}-u_{h\tau}\|_{\Lambda_h} \le C\|u_{h\tau}^{(k)}-u_{h\tau}\|_{L^2(\Gamma_1)} \le C\|u_h^{(k)} - u_h\|_{H^1(\Omega)^2}$
	in view of Lemmas \ref{Lem3.4}(i) and \ref{2.2}(i), we obtain
	\begin{equation}
		\|\lambda_h^{(k+1)}-\lambda_h\|_{\Lambda_h}^2 \le \|\lambda_h^{(k)}-\lambda_h\|_{\Lambda_h}^2 - (2\alpha\rho\ - C\rho^2)\|u_h^{(k)}-u_h\|_{H^1(\Omega)^2}^2, \label{3.37}
	\end{equation}
	and thus
	\begin{equation}
		(2\alpha\rho\ - C\rho^2)\|u_h^{(k)}-u_h\|_{H^1(\Omega)^2}^2 \le \|\lambda_h^{(k)}-\lambda_h\|_{\Lambda_h}^2 - \|\lambda_h^{(k+1)}-\lambda_h\|_{\Lambda_h}^2. \label{3.57}
	\end{equation}
	
	On the other hand, by virtue of Lemma \ref{Lem3.2}(i)(ii), we can choose $w_h \in V_{h,\sigma}$ such that $w_{h\tau}=\lambda_h^{(k)}-\lambda_h \text{ on }\Gamma_1$ and
	\begin{equation*}
		 \|w_h\|_{H^1(\Omega)^2} \le C\|\lambda_h^{(k)}-\lambda_h\|_{H^{1/2}(\Gamma_1)} \le C(h)\|\lambda_h^{(k)}-\lambda_h\|_{\Lambda_h},
	\end{equation*}
	where the constant $C(h)$ concerns the equivalence of the norms on the finite dimensional space $\Lambda_h$.
	Hence, it follows from (\ref{3.38}) with $v_h = w_h$ that
	\begin{align*}
		\|\lambda_h^{(k)}-\lambda_h\|_{\Lambda_h}^2 &= (w_{h\tau},\lambda_h^{(k)}-\lambda_h)_{\Lambda_h} = -a(\lambda_h^{(k)}-\lambda_h,w_h) \\
			&\le \|a\|\,\|u_h^{(k)}-u_h\|_{H^1(\Omega)^2}\|w_h\|_{H^1(\Omega)^2} \\
			&\le C(h)\|u_h^{(k)}-u_h\|_{H^1(\Omega)^2}\|\lambda_h^{(k)}-\lambda_h\|_{\Lambda_h},
	\end{align*}
	so that
	\begin{equation}
		-C(h)\|u_h^{(k)}-u_h\|_{H^1(\Omega)^2} \le -\|\lambda_h^{(k)}-\lambda_h\|_{\Lambda_h}. \label{3.39}
	\end{equation}
	Since the constant $C$ in (\ref{3.37}) is independent of $\rho$ (and even of $h$), if we choose $0 < \rho < \rho_0:=\frac{2\alpha}C$ then it follows from (\ref{3.37}) and (\ref{3.39}) that
	\begin{equation*}
		\|\lambda_h^{(k+1)}-\lambda_h\|_{\Lambda_h} \le \sqrt{1 - \frac{2\alpha\rho - C\rho^2}{C(h)^2}} \|\lambda_h^{(k)}-\lambda_h\|_{\Lambda_h},
	\end{equation*}
	where we may assume $C(h)^2 \ge 2\alpha\rho - C\rho^2$ (if not, take $C(h)^2 = \alpha^2/C$).
	Consequently, we conclude
	\begin{equation*}
		\|\lambda_h^{(k)}-\lambda_h\|_{\Lambda_h} \le \left(1 - \frac{2\alpha\rho - C\rho^2}{C(h)^2}\right)^\frac{k-1}2 \|\lambda_h^{(1)}-\lambda_h\|_{\Lambda_h} \to 0 \qquad (k\to\infty).
	\end{equation*}
	Then, from (\ref{3.57}) it also follows that $u_h^{(k)} \to u_h$ in $H^1(\Omega)^2$ as $k \to \infty$.
	
	Finally, subtracting (\ref{3.34}) from (\ref{3.13}) with test functions in $\stackrel{\circ\,}{V_h}$ gives
	\begin{equation*}
		a(u_h-u_h^{(k)},v_h) + b(v_h, p_h - p_h^{(k)}) = 0 \qquad (\forall v_h \in \stackrel{\circ\,}{V_h}).
	\end{equation*}
	Therefore, by Lemma \ref{Lem3.1}(i) we have
	\begin{align*}
		\beta\|p_h^{(k)}-p_h\|_{L^2(\Omega)} &\le \sup_{v_h\in \stackrel{\circ\,}{V_h}}\frac{b(v_h,p_h^{(k)}-p_h)}{\|v_h\|_{H^1(\Omega)^2}} = \sup_{v_h\in \stackrel{\circ\,}{V_h}}\frac{-a(u_h^{(k)}-u_h,v_h)}{\|v_h\|_{H^1(\Omega)^2}} \\
			&\le \|a\|\,\|u_h^{(k)}-u_h\|_{H^1(\Omega)^2} \\
			&\to 0 \qquad (k \to \infty).
	\end{align*}
	This completes the proof. 
\end{proof}

\section{Discretization of the Stokes problem with LBCF}\label{Sec5}
\subsection{Existence and uniqueness results}\label{Sec5.2}
Approximate problems to Problem VI (therefore, Problem PDE) in the case of LBCF are as follows.
\begin{Prob}[VI$_h$]
	Find $(u_h, p_h) \in V_{\tau h}\times Q_h$ such that
	\begin{numcases}{\hspace{-7mm}}
		\hspace{-2mm} a(u_h,v_h\!-\!u_h) \!+\! b(v_h\!-\!u_h,p_h) \!+\! j_h(\!v_{hn}\!) \!-\! j_h(\!u_{hn}\!) \!\ge\! (f,v_h\!-\!u_h)\!_{L^2(\Omega)^2} &\hspace{-7mm} $(\forall v_h \!\!\in\!\! V_{\tau h})$\!, \label{4.23} \\
		\hspace{-2mm} b(u_h,q_h) = 0 &\hspace{-7mm} $(\forall q_h \!\!\in\!\! Q_h)$. \label{4.24}
	\end{numcases}
\end{Prob}
\begin{Prob}[VI$_{h,\sigma}$]
	Find $u_h\in V_{\tau h,\sigma}$ such that
	\begin{equation}
		a(u_h,v_h-u_h) + j_h(v_{hn}) - j_h(u_{hn}) \ge (f,v_h-u_h)_{L^2(\Omega)^2} \qquad (\forall v_h\in V_{\tau h,\sigma}). \label{4.8}
	\end{equation}
\end{Prob}
\begin{Prob}[VE$_h$]
	Find $(u_h, p_h, \lambda_h) \in V_{\tau h}\times Q_h\times\tilde\Lambda_h$ such that
	\begin{numcases}{\hspace{-1cm}}
		a(u_h,v_h) + b(v_h,p_h) + (v_{hn}, \lambda_h)_{\Lambda_h} = (f,v_h)_{L^2(\Omega)^2} & $(\forall v_h\in V_{\tau h})$, \label{4.13} \\
		b(u_h,q_h) = 0 & $(\forall q_h\in Q_h)$, \label{4.14} \\
		(u_{hn}, \mu_h - \lambda_h)_{\Lambda_h} \le 0 & $(\forall\mu_h \in \tilde\Lambda_h)$. \label{4.15}
	\end{numcases}
\end{Prob}

\begin{Thm}\label{Thm4.1}
	{\rm(i)} Problem VI$_{h,\sigma}$ admits a unique solution $u_h\in V_{\tau h,\sigma}$.
	Furthermore, it satisfies the following equation:
	\begin{equation}
		a(u_h,u_h) + j_h(u_{hn}) = (f,u_h)_{L^2(\Omega)^2}. \label{4.12}
	\end{equation}
	
	\noindent{\rm(ii)} Problems $\mathrm{VI}_h$ and $\mathrm{VE}_h$ are equivalent in the following sense.
	
	{\rm(a)} If $(u_h,p_h)\in V_{\tau h}\times Q_h$ is a solution of Problem {\rm VI}$_h$, then there exists a unique $\lambda_h\in \tilde\Lambda_h$
		such that $(u_h,p_h,\lambda_h)$ solves Problem {\rm VE}$_{h}$.
		
	{\rm(b)} If $(u_h,p_h,\lambda_h) \in V_{\tau h}\times Q_h\times\tilde\Lambda_h$ is a solution of Problem {\rm VE}$_h$, then $(u_h,p_h)$ solves Problem {\rm VI}$_h$.
\end{Thm}
\begin{proof}
	Statements (i) and (ii)(a) are proved by the same way as Theorems \ref{Thm3.1}(i) and (ii)(b), respectively.
	We obtain (ii)(b) by a discussion similar to that in the proof of Theorem \ref{Thm3.1}(ii)(c). 
\end{proof}
\begin{Rem}\label{Rem4.2}
	It is clear that if $(u_h,p_h)$ is a solution of Problem VI$_h$, then $u_h$ is that of Problem VI$_{h,\sigma}$.
	However, the "converse" is no longer true. In fact, the pressure, especially its additive constant, need not to be uniquely determined
	even if the solution $u_h$ of Problem VI$_{h,\sigma}$ is given. This situation is quite different from the case of SBCF.
\end{Rem}

The next theorem guarantees the existence of a solution of Problem VI$_h$.
\begin{Thm}\label{Thm4.2}
	{\rm(i)} There exists at least one solution of Problem {\rm VE}$_h$, and the velocity part is unique.
	
	{\rm(ii)} If $(u_h, p_h, \lambda_h)$ and $(u_h, p_h^*, \lambda_h^*)$ are two solutions of Problem {\rm VE}$_h$, then there exists a unique $\delta_h\in\mathbf R$ such that
	\begin{equation}
		p_h = p_h^* + \delta_h \qquad{\rm and}\qquad \lambda_h(M) = \lambda_h^*(M) + \frac{\delta_h}{g(M)} \quad (\forall M \!\in\, \stackrel\circ\Gamma_{1,h}). \label{4.31}
	\end{equation}
	
	{\rm(iii)} Under the assumptions in {\rm(ii)}, if we suppose $u_{hn}\neq 0$ on $\Gamma_1$, then $\delta_h=0$. Namely, a solution of Problem $\mathrm{VE}_h$ is unique.
\end{Thm}
\begin{proof}
	(i) The uniqueness of the velocity is obvious; see Remark \ref{Rem4.2}
	
	To prove the existence, let $u_h \in V_{h,\sigma}$ be the solution of {\rm VI}$_{h,\sigma}$.
	Taking $u_h \pm v_h$ as a test function in (\ref{4.8}), with an arbitrary $v_h \in \stackrel{\circ\,}{V_h}\cap\; V_{\tau h,\sigma}$, we obtain
	\begin{equation*}
		a(u_h,v_h) = (f,v_h)_{L^2(\Omega)^2} \qquad (\forall v_h \in \stackrel{\circ\,}{V_h}\cap\; V_{\tau h,\sigma}).
	\end{equation*}
	Hence, from Lemma \ref{Lem3.1}(i), we deduce the unique existence of $p_h^0 \in\; \stackrel{\circ~}{Q_h}$ such that
	\begin{equation*}
		a(u_h,v_h) + b(v_h,p_h^0) = (f,v_h)_{L^2(\Omega)^2} \qquad (\forall v_h \in \stackrel{\circ\,}{V_h}).
	\end{equation*}
	Moreover, noting that $\{v_h\in V_{\tau h} \,|\, (v_{hn}, \eta_h)_{\Lambda_h}=0\, (\forall\eta_h\in\Lambda_h)\} = \stackrel{\circ\,}{V_h}$, it follows from
	Lemma \ref{Lem3.2}(iii) that there exists a unique $\lambda_h^0 \in \Lambda_h$ satisfying
	\begin{equation}
		a(u_h,v_h) + b(v_h,p_h^0) + (v_{hn}, \lambda_h^0)_{\Lambda_h} = (f,v_h)_{L^2(\Omega)^2} \qquad (\forall v_h\in V_{\tau h}). \label{4.9}
	\end{equation}
	
	For every $\eta_h\in \Lambda_h\cap L^2_0(\Gamma_1)$, by Lemma \ref{Lem3.2}(ii) we can choose $w_h \in V_{\tau h,\sigma}$ such that
	\begin{equation*}
		w_{hn} = \eta_h \quad\text{on}\quad \Gamma_1.
	\end{equation*}
	Hence, taking $v_h = w_h$ in (\ref{4.8}) and (\ref{4.9}), we obtain
	\begin{gather}
		a(u_h, w_h-u_h) + j_h(\eta_h) - j_h(u_{hn}) \ge (f, w_h-u_h)_{L^2(\Omega)^2}, \label{4.10} \\
		a(u_h, w_h) + (\eta_h, \lambda_h^0)_{\Lambda_h} = (f, w_h)_{L^2(\Omega)^2}. \label{4.11}
	\end{gather}
	Substituting (\ref{4.12}) and (\ref{4.11}) into (\ref{4.10}) gives
	\begin{equation*}
		(\eta_h, \lambda_h^0)_{\Lambda_h} \le j_h(\eta_h) \qquad (\forall \eta_h\in \Lambda_h\cap L^2_0(\Gamma_1)).
	\end{equation*}
	We apply Hahn-Banach's theorem to deduce the existence of some $\lambda_h \in \Lambda_h$ such that
	\begin{equation*}
		(\eta_h, \lambda_h)_{\Lambda_h} \le j_h(\eta_h) \qquad (\forall \eta_h\in \Lambda_h).
	\end{equation*}
	Therefore, Lemma \ref{Lem3.3}(ii) implies that $\lambda_h \in \tilde\Lambda_h$.
	Furthermore, since $\lambda_h$ satisfies
	\begin{equation*}
		(\eta_h, \lambda_h - \lambda_h^0)_{\Lambda_h} = 0 \qquad (\forall \eta_h\in \Lambda_h\cap L^2_0(\Gamma_1)),
	\end{equation*}
	it follows from Lemma \ref{Lem3.3}(iv) that there exists some $\delta_h \in \mathbf R$ such that
	\begin{equation*}
		\lambda_h(M) = \lambda_h^0(M) + \frac{\delta_h}{g(M)} \qquad (\forall M \!\in\, \stackrel\circ\Gamma_{1,h}).
	\end{equation*}
	Thus, from Simpson's formula and (\ref{4.9}), we obtain
	\begin{align}
		&(f,v_h)_{L^2(\Omega)^2} \notag \\
			= &a(u_h,v_h) + b(v_h,p_h^0) + (v_{hn}, \lambda_h^0)_{\Lambda_h} \label{4.17} \\
			= &a(u_h,v_h) + b(v_h,p_h^0) + (v_{hn}, \lambda_h)_{\Lambda_h} - \delta_h\sum_{i=1}^m\frac{|e_i|}6 \left( v_{hn,i} + 4v_{hn,i+\frac12} + v_{hn,i+1} \right) \\
			= &a(u_h,v_h) + b(v_h,p_h^0) + (v_{hn}, \lambda_h)_{\Lambda_h} - \delta_h \int_{\Gamma_1}\!v_{hn}\,ds \\
			= &a(u_h,v_h) + b(v_h,p_h^0 + \delta_h) + (v_{hn}, \lambda_h)_{\Lambda_h} \quad (\forall v_h \in V_{\tau h}). \label{4.16}
	\end{align}
	This establishes (\ref{4.13}) if we define $p_h=p_h^0+\delta_h$.
	
	Equation (\ref{4.14}) obviously holds because $u_h\in V_{h,\sigma}$.
	It remains to show (\ref{4.15}), which is equivalent to $(u_{hn}, \lambda_h)_{\Lambda_h} = j_h(u_{hn})$ by Lemma \ref{Lem3.3}(ii).
	This is indeed obtained from (\ref{4.16}) with $v_h=u_h$ and (\ref{4.12}).
	
	(ii) Let $(u_h, p_h, \lambda_h)$ and $(u_h, p_h^*, \lambda_h^*)$ be two solutions of Problem VE$_h$.
	Because the uniqueness of the pressure up to additive constants is shown in the proof of (i), there exists a unique constant $\delta_h$ such that $p_h=p_h^*+\delta_h$.

	Since $(u_h, p_h, \lambda_h)$ and $(u_h, p_h^*, \lambda_h^*)$ satisfy (\ref{4.13}), subtracting the two equations and calculating
	in a manner similar to (\ref{4.17})--(\ref{4.16}), we obtain
	\begin{equation}
		(v_{hn}, \hat \lambda_h)_{\Lambda_h} = 0 \qquad (\forall v_h\in V_{h\tau}). \label{4.30}
	\end{equation}
	Here, $\hat \lambda_h\in\Lambda_h$ is defined by $\hat \lambda_h(M) = \lambda_h(M) - \lambda_h^*(M) - \delta_h/g(M)$ for $M\in\; \stackrel\circ\Gamma_{1,h}$.
	It follows from (\ref{4.30}) together with Lemma \ref{Lem3.2}(i) that $(\eta_h,\hat \lambda_h)_{\Lambda_h}=0$ for all $\eta_h\in\Lambda_h$. 
	Hence $\hat\lambda_h=0$, and (\ref{4.31}) is proved.
	
	(iii) The assumption $u_{hn} \neq 0$ implies that either of the following is true:
	
	\quad(a) There exists $M \!\in\, \stackrel\circ\Gamma_{1,h}$ such that $u_{hn}(M) > 0$.
	
	\quad(b) There exists $M \!\in\, \stackrel\circ\Gamma_{1,h}$ such that $u_{hn}(M) < 0$.
	
	\noindent Since $u_h\in V_{\tau h,\sigma}$, we see that
	\begin{equation*}
		\int_{\Gamma_1}u_{hn}\,ds = \int_{\Omega}\mathrm{div}\,u_h\,dx = 0.
	\end{equation*}
	Therefore, both (a) and (b) above must hold true.
	Then, from (\ref{4.15}) combined with Lemma \ref{Lem3.3}(ii), we find some $M,N \in\; \stackrel\circ\Gamma_{1,h}$ satisfying $\lambda_h(M)=1$ and $\lambda_h(N)=-1$.
	Consequently, the additive constant $\delta_h$ appearing in (ii) cannot attain any value except $0$ because $\lambda_h \in \tilde\Lambda_h$.
	This completes the proof of Theorem \ref{Thm4.2}.
\end{proof}

\subsection{Error analysis}\label{Sec5.3}
Let us begin with the following analogue of Proposition \ref{Prop3.1}.
\begin{Prop}\label{Prop4.1}
	Let $(u,p)$ be a solution of Problem {\rm VI} and $(u_h,p_h)$ be that of Problem {\rm VI}$_h$\, for $h>0$. Then,
	
	{\rm(i)} it holds that
	\begin{equation*}
		\|u_h\|_{H^1(\Omega)^2} \le \|f\|_{L^2(\Omega)^2}/\alpha.
	\end{equation*}
	
	{\rm(ii)} for every $v_h\in V_{\tau h}$ and $q_h\in Q_h$, it holds that
	\begin{align}
		&\alpha\|u-u_h\|_{H^1(\Omega)^2}^2 \le a(u-u_h,u-v_h) + b(u_h-u,p-q_h) + b(v_h-u,p_h-p) \notag \\
		&\hspace{1.5cm} + (\sigma_n,v_{hn}-u_n)_{L^2(\Gamma_1)} + j(u_{hn}) - j_h(u_{hn}) + j_h(v_{hn}) - j(u_n). \label{4.21}
	\end{align}
	
	{\rm(iii)} for every $q_h \in\, \stackrel{\circ~}{Q_h}$, it holds that
	\begin{align}
		\|p^0-p_h^0\|_{L^2(\Omega)} \le \left( 1 + \frac{\|b\|}\beta \right) \|p^0-q_h\|_{L^2(\Omega)^2} + \frac{\|a\|}\beta \|u-u_h\|_{H^1(\Omega)^2}, \label{4.20}
	\end{align}
	where $p^0=p-(p,1)_{L^2(\Omega)}/|\Omega|$ \;and\; $p_h^0=p_h-(p_h,1)_{L^2(\Omega)}/|\Omega|$.
	
	{\rm(iv)} for every $q_h\in Q_h$, it holds that
	\begin{align}
		\|p-p_h\|_{L^2(\Omega)} \le \left( 1 + \frac{\|b\|}\beta \right) \|p-q_h\|_{L^2(\Omega)} + \frac{\|a\|}\beta \|u-u_h\|_{H^1(\Omega)^2} + C. \label{4.19}
	\end{align}
\end{Prop}
\begin{proof}
	Statements other than (iv) can be proved by the same way as Proposition \ref{Prop3.1}.
	To show (iv), we let $q_h\in Q_h$. It is clear that $\|p-p_h\|_{L^2(\Omega)} \le \|p-q_h\|_{L^2(\Omega)} + \|q_h-p_h\|_{L^2(\Omega)}$.
	To bound the latter term, we deduce from Lemma \ref{Lem3.1}, together with (\ref{2.13}) and (\ref{4.13}), that
	\begin{align}
		&\beta\|q_h-p_h\|_{L^2(\Omega)} \le \sup_{v_h\in V_{\tau h}}\frac{b(v_h, q_h-p_h)}{\|v_h\|_{H^1(\Omega)^2}} = \sup_{v_h\in V_{\tau h}}\frac{b(v_h, q_h-p) + b(v_h, p-p_h)}{\|v_h\|_{H^1(\Omega)^2}} \notag \\
		&\hspace{1cm}	= \sup_{v_h\in V_{\tau h}}\frac{b(v_h, q_h-p) - a(u-u_h,v_h) + (\sigma_n, v_{hn})_{L^2(\Gamma_1)} + (v_{hn}, \lambda_h)_{\Lambda_h}}{\|v_h\|_{H^1(\Omega)^2}} \notag \\
		&\hspace{1cm}	\le \|b\|\,\|p - q_h\|_{L^2(\Omega)} + \|a\|\,\|u - u_h\|_{H^1(\Omega)^2} + C. \label{4.18}
	\end{align}
	Here, to derive (\ref{4.18}), we have used the estimates
	\begin{gather*}
		|(\sigma_n, v_{hn})_{L^2(\Gamma_1)}| \le \|\sigma_n\|_{L^2(\Gamma_1)} \|v_{hn}\|_{L^2(\Gamma_1)} \le C\|v_{hn}\|_{H^\frac12(\Gamma_1)} \le C\|v_h\|_{H^1(\Omega)^2}, \\
		|(v_{hn}, \lambda_h)_{\Lambda_h}| \le j_h(v_{hn}) \le C\|v_{hn}\|_{L^2(\Gamma_1)} \le C\|v_{hn}\|_{H^\frac12(\Gamma_1)} \le C\|v_h\|_{H^1(\Omega)^2},
	\end{gather*}
	which are obtained from Lemmas \ref{Lem2.2}(i), \ref{Lem3.3}(i), and \ref{Lem3.4}(i).
	The desired inequality (\ref{4.19}) immediately follows from (\ref{4.18}).
\end{proof}
\begin{Rem}
	If $p-p_h \in L^2_0(\Omega)$ and $q_h\in\; \stackrel{\circ~}{Q_h}$, we can take $C=0$ in (\ref{4.19}) according to the equality $p-p_h = p^0-p_h^0$ combined with (\ref{4.20}).
\end{Rem}

We state the rate of convergence result for the case of LBCF, which is not better than that of SBCF because of the influence of an additive constant of the pressure.
\begin{Thm}\label{Thm4.3}
	Let $(u,p)$ be a solution of Problem {\rm VI} and $(u_h,p_h)$ be that of Problem {\rm VI}$_h$\, for $0\!<\!h\!<\!1$,
	and suppose $(u,p) \in H^{1+\epsilon}(\Omega)^2 \times H^\epsilon(\Omega)$ with $0< \epsilon\le 2$. Then,
	\begin{equation}
		\|u - u_h\|_{H^1(\Omega)^2} + \|p^0 - p_h^0\|_{L^2(\Omega)} \le Ch^{\min\{\frac\epsilon2, \frac14\}},
	\end{equation}
	where $p^0=p-(p,1)_{L^2(\Omega)}/|\Omega|$ \,and\, $p_h^0=p_h-(p_h,1)_{L^2(\Omega)}/|\Omega|$.
\end{Thm}
\begin{proof}
	Let us take $v_h=\mathcal I_hu$ and $q_h=\Pi_hp$ in (\ref{4.21}) and bound from above each term on the right-hand side.
	By (\ref{4.19}), we have
	\begin{align}
		|b(\mathcal I_hu-u, p_h-p)| &\le \|b\|\,\|\mathcal I_hu-u\|_{H^1(\Omega)^2}\|p_h-p\|_{L^2(\Omega)} \\
			&\le C( h^\epsilon\|u-u_h\|_{H^1(\Omega)^2} + h^{2\epsilon} + h^\epsilon ).
	\end{align}
	For the other terms, we employ the same estimates as those in the proof of Theorem \ref{Thm3.2}.
	Then it follows that
	\begin{equation}
		\|u - u_h\|_{H^1(\Omega)^2}^2 \le C( h^\epsilon\|u-u_h\|_{H^1(\Omega)^2} + h^{2\epsilon} + h^\epsilon + h^\frac12 ),
	\end{equation}
	which implies
	\begin{equation}
		\|u - u_h\|_{H^1(\Omega)^2} \le Ch^{\min\{\frac\epsilon2, \frac14\}}. \label{4.22}
	\end{equation}
	$\|p^0 - p_h^0\|_{L^2(\Omega)} \le Ch^{\min\{\frac\epsilon2, \frac14\}}$ follows from (\ref{4.20}) and (\ref{4.22}) and this completes the proof.
\end{proof}
\begin{Rem}
	(i) If we assume, in addition, that $p-p_h \in L^2_0(\Omega)$ then we can establish the result of $O(h^{\min\{\epsilon, \frac14\}})$.
	Moreover, as in the case of SBCF, under suitable conditions regarding the signs of $u_{hn}$ and $(\mathcal I_hu)_n$ on $\Gamma_1$, it can be improved to $O(h^{\min\{\epsilon, 1\}})$,
	or even $O(h^\epsilon)$ if $g$ is affine.
	
	(ii) When the uniqueness of Problem VI holds, we can obtain a strong convergence result for the error of the pressure including the additive constant.
	In fact, the uniform boundedness of $p_h$ in $L^2(\Omega)$ gives a weak convergence limit $\tilde p$ for some subsequence $p_{h'}$.
	Since $u_h\to u$ in $H^1(\Omega)^2$, we have
	\begin{equation*}
		|j_h(u_{hn}) - j(u_n)| \le |j_h(u_{hn}) - j(u_{hn})| + |j(u_{hn}) - j(u_n)| \to 0 \quad (h\to 0).
	\end{equation*}
	Therefore, taking the limit $h'\to 0$ in (\ref{4.23})--(\ref{4.24}), we find that $(u,\tilde p)$ is a solution of Problem VI.
	Hence $p=\tilde p$, and from $p_{h'}\rightharpoonup p$ and $p_h^0\to p^0$, we conclude the strong convergence of the whole sequence.
\end{Rem}

\subsection{Numerical realization}\label{Sec5.4}
Based on Problem VE$_h$, we propose the following Uzawa-type method to compute the approximate solution $(u_h,p_h)$ numerically.
\begin{Alg}\label{Alg4.1}
	Choose an arbitrary $\lambda_h^{(1)} \in \tilde\Lambda_h$ and $\rho > 0$. Iterate the following two steps for $k=1,2,\cdots:$
	
	{\rm\bf \underline{Step 1}}\; With $\lambda_h^{(k)}$ known, determine $(u_h^{(k)},p_h^{(k)}) \in V_{\tau h}\times Q_h$ by
	\begin{numcases}{\hspace{-5mm}}
		\!\!a(u_h^{(k)},v_h) + b(v_h,p_h^{(k)}) = (f,v_h)_{L^2(\Omega)^2} - (v_{hn}, \lambda_h^{(k)})_{\Lambda_h} &\hspace{-2mm}$(\forall v_h \!\in\! V_{\tau h})$, \label{4.25} \\
		\!\!b(u_h^{(k)},q_h) = 0 &$\hspace{-2mm}(\forall q_h\!\in\! Q_h)$. \label{4.32}
	\end{numcases}
	
	{\rm\bf \underline{Step 2}}\; Renew $\lambda_h^{(k+1)} \in \tilde\Lambda_h$ by
	\begin{equation}
		\lambda_h^{(k+1)} = \mathrm{Proj}_{\tilde\Lambda_h}(\lambda_h^{(k)} + \rho u_{hn}^{(k)}). \label{4.26}
	\end{equation}
\end{Alg}
\begin{Rem}
	The unique existence of $(u_h^{(k)},p_h^{(k)})$ satisfying (\ref{4.25}) and (\ref{4.32}) is guaranteed by the inf-sup condition in Lemma \ref{Lem4.2}.	
\end{Rem}

\begin{Thm}\label{Thm4.4}
	Under the same notation as {\rm Algorithm \ref{Alg4.1}}, there exists a constant $\rho_0>0$ independent of $h$ such that if $\rho$ satisfies $0<\rho<\rho_0$, then
	$(u_h^{(k)},p_h^{(k)},\lambda_h^{(k)})$ converges to some solution of Problem $\mathrm{VE}_h$ in $H^1(\Omega)^2\times L^2(\Omega)\times \Lambda_h$, as $k \to \infty$.
\end{Thm}
\begin{proof}
	First we show the boundedness of the sequence $\{(u_h^{(k)},p_h^{(k)},\lambda_h^{(k)})\}_k$.
	In fact, taking $v_h=u_h$ in (\ref{4.25}), we find from Korn's inequality (\ref{2.5}) that
	\begin{align*}
		\alpha\|u_h^{(k)}\|_{H^1(\Omega)^2}^2 &\le a(u_h^{(k)},u_h^{(k)}) = (f,u_h^{(k)})_{L^2(\Omega)^2} - (u_{hn}^{(k)},\lambda_h^{(k)})_{\Lambda_h} \\
			&\le \|f\|_{L^2(\Omega)^2}\|u_h^{(k)}\|_{L^2(\Omega)^2} + |(u_{hn}^{(k)},\lambda_h^{(k)})_{\Lambda_h}| \\
			&\le C\|u_h^{(k)}\|_{H^1(\Omega)^2},
	\end{align*}
	which gives $\|u_h^{(k)}\|_{H^1(\Omega)^2} \le C$. Here, to derive the last line, we have used
	\begin{equation}
		|(u_{hn}^{(k)},\lambda_h^{(k)})_{\Lambda_h}| \le j_h(u_{hn}^{(k)}) \le C\|u_{hn}^{(k)}\|_{L^2(\Gamma_1)} \le C\|u_h^{(k)}\|_{H^1(\Omega)^2}, \label{4.33}
	\end{equation}
	which is obtained from Lemmas \ref{Lem3.3}(i), \ref{Lem3.4}(i), and \ref{Lem2.2}(i).
	Then Lemma \ref{Lem4.2}, together with (\ref{4.25}), implies that
	\begin{align*}
		\beta\|p_h^{(k)}\|_{L^2(\Omega)^2} &\le \sup_{v_h\in V_{\tau h}}\frac{b(v_h,p_h^{(k)})}{\|v_h\|_{H^1(\Omega)^2}} \\
		&= \sup_{v_h\in V_{\tau h}}\frac{(f,v_h)_{L^2(\Omega)^2} - a(u_h^{(k)},v_h) - (v_{hn},\lambda_h^{(k)})_{\Lambda_h}}{\|v_h\|_{H^1(\Omega)^2}} \\
		&\le \|f\|_{L^2(\Omega)^2} + \|a\|\,\|u_h^{(k)}\|_{H^1(\Omega)^2} + \sup_{v_h\in V_{\tau h}} \frac{|(v_{hn},\lambda_h^{(k)})_{\Lambda_h}|}{\|v_h\|_{H^1(\Omega)^2}} \\
		&\le C,
	\end{align*}
	where $|(v_{hn},\lambda_h^{(k)})_{\Lambda_h}|$ in the third line is estimated in a manner similar to (\ref{4.33}).
	It is clear that $\{ \lambda_h^{(k)} \}_k$ is bounded because $\lambda_h^{(k)}\in \tilde\Lambda_h$.
	
	\hspace{-1mm}Therefore, we can extract a subsequence $\{(u_h^{(k')},p_h^{(k')},\lambda_h^{(k')})\}_{k'}$ converging to some element $(u_h,p_h,\lambda_h)\in V_{\tau h}\times Q_h\times\tilde\Lambda_h$.
	Making $k=k'$ and $k'\to\infty$ in (\ref{4.25})--(\ref{4.26}), we obtain
	\begin{numcases}{\hspace{-1cm}}
		a(u_h,v_h) + b(v_h,p_h) + (v_{hn}, \lambda_h)_{\Lambda_h} = (f,v_h)_{L^2(\Omega)^2} & $(\forall v_h\in V_{\tau h})$, \notag \\
		b(u_h,q_h) = 0 & $(\forall q_h\in Q_h)$, \notag \\
		\lambda_h = \mathrm{Proj}_{\tilde\Lambda_h}(\lambda_h + \rho u_{hn}). \notag
	\end{numcases}
	Consequently, since the last equation is equivalent to (\ref{4.15}) by virtue of Lemma \ref{Lem3.3}(ii), we see that $(u_h,p_h,\lambda_h)$ is a solution of Problem VE$_h$.
	
	It remains only to prove that the whole sequence converges to $(u_h,p_h,\lambda_h)$.
	Subtracting (\ref{4.25}) from (\ref{4.13}), we have
	\begin{equation}
		a(u_h-u_h^{(k)}, v_h) + b(v_h, p_h-p_h^{(k)}) + (v_{hn}, \lambda_h-\lambda_h^{(k)})_{\Lambda_h} = 0 \quad (\forall v_h\in V_{\tau h}). \label{4.28}
	\end{equation}
	In particular, if we take $v_h = u_h - u_h^{(k)}$, then
	\begin{equation*}
		(u_{hn} - u_{hn}^{(k)}, \lambda_h - \lambda_h^{(k)})_{\Lambda_h} = -a(u_h-u_h^{(k)}, u_h-u_h^{(k)}).
	\end{equation*}
	Therefore, it follows from a general property of the projection operator $\mathrm{Proj}_{\Lambda_h}$ that
	\begin{align}
		\|\lambda_h - \lambda_h^{(k+1)}\|_{\Lambda_h}^2 &\le \|\lambda_h - \lambda_h^{(k)} + \rho (u_{hn} - u_{hn}^{(k)})\|_{\Lambda_h}^2 \notag \\
			&= \|\lambda_h - \lambda_h^{(k)}\|_{\Lambda_h}^2 \!- 2\rho a(u_{hn} - u_{hn}^{(k)}, u_{hn} - u_{hn}^{(k)}) + \rho^2\|u_{hn} - u_{hn}^{(k)}\|_{\Lambda_h}^2 \notag \\
			&\le \|\lambda_h - \lambda_h^{(k)}\|_{\Lambda_h}^2 + (C\rho^2 - 2\alpha\rho)\|u_h - u_h^{(k)}\|_{H^1(\Omega)^2}^2, \label{4.29}
	\end{align}
	where we have applied Korn's inequality (\ref{2.5}) and the estimate $\|u_{hn} - u_{hn}^{(k)}\|_{\Lambda_h} \le C\|u_{hn} - u_{hn}^{(k)}\|_{L^2(\Gamma_1)} \le C\|u_h - u_h^{(k)}\|_{H^1(\Omega)^2}$
	to derive the last line.
	Since the constant $C$ in (\ref{4.29}) is independent of $\rho$ (and even of $h$), we choose $0<\rho<\rho_0:=\frac{2\alpha}{C}$ to obtain
	\begin{equation}
		0 \le \|u_h - u_h^{(k)}\|_{H^1(\Omega)^2}^2 \le \frac1{2\alpha\rho - C\rho^2}\bigg( \|\lambda_h - \lambda_h^{(k)}\|_{\Lambda_h}^2 - \|\lambda_h - \lambda_h^{(k+1)}\|_{\Lambda_h}^2 \bigg). \label{4.27}
	\end{equation}
	Hence the sequence $\{ \|\lambda_h - \lambda_h^{(k)}\|_{\Lambda_h} \}_k$ is decreasing.
	Noting that a decreasing sequence in $\mathbf R$, bounded from below, converges to its infimum and that $\inf_{k}\|\lambda_h - \lambda_h^{(k)}\|_{\Lambda_h}=0$ by construction,
	we conclude that $\lambda_h^{(k)}\to \lambda_h$ in $\Lambda_h$ as $k\to\infty$. From this and (\ref{4.27}), it also follows that $u_h^{(k)}\to u_h$ in $H^1(\Omega)^2$.
	Finally, from the inf-sup condition given in Lemma \ref{Lem4.2} combined with (\ref{4.28}), we have
	\begin{align*}
		\beta\|p_h-p_h^{(k)}\|_{L^2(\Omega)} &\le \sup_{v_h\in V_{\tau h}}\frac{b(v_h,p_h-p_h^{(k)})}{\|v_h\|_{H^1(\Omega)^2}} \\
		&= \sup_{v_h\in V_{\tau h}}\frac{-a(u_h-u_h^{(k)},v_h) - (v_{hn},\lambda_h-\lambda_h^{(k)})_{\Lambda_h}}{\|v_h\|_{H^1(\Omega)^2}} \\
		&\le \|a\|\;\|u_h^{(k)}-u_h\|_{H^1(\Omega)^2} + \|\lambda_h^{(k)}-\lambda_h\|_{\Lambda_h} \sup_{v_h\in V_{\tau h}}\frac{\|v_{hn}\|_{\Lambda_h}}{\|v_h\|_{H^1(\Omega)^2}} \\
		&\le \|a\|\;\|u_h^{(k)}-u_h\|_{H^1(\Omega)^2} + C\|\lambda_h^{(k)}-\lambda_h\|_{\Lambda_h} \to 0 \qquad (k \to \infty).
	\end{align*}
	This completes the proof. 
\end{proof}
\begin{Rem}
	(i) The resulting solution of Problem VE$_h$ as the limit of $(\!u_h^{(k)}\!, p_h^{(k)}\!, \lambda_h^{(k)}\!)$, especially its additive constant of the pressure,
	may depend on a choice of $0<\rho<\rho_0$ or that of the starting value $\lambda_h^{(1)}$.
	However, if $u_{hn} \neq 0$ on $\Gamma_1$, and hence the uniqueness of the solution of Problem VE$_h$ is valid, then it is obviously independent of them.

	(ii) Contrary to the case of SBCF, it is difficult to prove an exponential convergence of the iterative solution $(u_h^{(k)},p_h^{(k)},\lambda_h^{(k)})$ because
	we do not know whether $\lambda_h - \lambda_h^{(k)}\in L^2_0(\Gamma_1)$, which is necessary to deduce an extension of $\lambda_h - \lambda_h^{(k)}$ to $V_{\tau h,\sigma}$.
\end{Rem}

\section{Numerical examples}\label{Sec6}
We assume $\Omega = (0,1)^2$, the boundary of which consists of two portions $\Gamma_0$ and $\Gamma_1$ given by
\begin{gather}
	\Gamma_0 = \{(0,y)\,|\,0<y<1\} \cup \{(x,0)\,|\,0\le x\le 1\} \cup \{(1,y)\,|\,0<y<1\}, \\
	\Gamma_1 = \{(x,1)\,|\,0<x<1\}.
\end{gather}
In particular, the set of extreme points is $\overline\Gamma_0\cap\overline\Gamma_1 = \{(0,1),(1,1)\}$.
For the triangulation $\mathscr T_h$ of $\overline\Omega$, we employ a uniform $N\times N$ Friedrichs-Keller type mesh, where $N$ denotes the division number of each side of the square $\overline\Omega$.

Let us consider
\begin{equation}
	\begin{cases}\label{5.1}
		u_1(x,y) &= 20x^2(1-x)^2y(1-y)(1-2y), \\
		u_2(x,y) &= -20x(1-x)(1-2x)y^2(1-y)^2, \\
		p(x,y) &= 40x(1-x)(1-2x)y(1-y)(1-2y) \\
		&\hspace{2cm}	+ 4(6x^5-15x^4+10x^3)(2y-1) - 2,
	\end{cases}
\end{equation}
which turns out to be the solution of the Stokes equations under the adhesive boundary condition for $\nu=1$ and $f$ given by
\begin{numcases}{\hspace{-2mm}}
	\hspace{-2mm} f_1(x,y) \!=\! 0, \notag \\
	\hspace{-2mm} f_2(x,y) \!=\! 120(\!2x\!-\!1\!)y^2(\!1\!-\!y\!)^2 \!+\! 80x(\!1\!-\!x\!)(\!1\!-\!2x\!)(\!6y^2\!-\!6y\!+\!1\!) \!+\! 8(\!6x^5\!-\!15x^4\!+\!10x^3\!). \notag
\end{numcases}
By direct computation, we have
\begin{gather}
	\max_{\overline\Gamma_1} |\sigma_\tau| = \max_{0\le x\le 1} \left|20x^2(1-x)^2\right| = \frac54 = 1.25, \\
	\max_{\overline\Gamma_1} |\sigma_n| = \max_{0\le x\le 1} \left|-4(6x^5-15x^4+10x^3) + 2\right| = 2.
\end{gather}
Now, if we impose SBCF or LBCF on $\Gamma_1$, with $g$ being constant, instead of the adhesive boundary condition, then in the case of SBCF, we find that
\begin{equation*}
	\begin{cases}
		g \ge 1.25 &\Longrightarrow \text{ (\ref{5.1}) remains a solution.} \\
		g < 1.25 &\Longrightarrow \text{ (\ref{5.1}) is no longer a solution and a non-trivial slip occurs.}
	\end{cases}
\end{equation*}
and in the case of LBCF,
\begin{equation*}
	\begin{cases}
		g \ge 2 &\Longrightarrow \text{ (\ref{5.1}) remains a solution.} \\
		g < 2 &\Longrightarrow \text{ (\ref{5.1}) is no longer a solution and a non-trivial leak occurs.}
	\end{cases}
\end{equation*}
We indeed observe some of the abovementioned phenomena in our numerical computation, as indicated in the plots of the velocity field shown below in Figures \ref{Fig5.1} and \ref{Fig5.2}.
In addition, we find that the bigger (resp. smaller) the threshold $g$ of a tangential or normal stress becomes, the more difficult (resp. easier) it becomes
for a non-trivial slip or leak to occur, which is in agreement with our natural intuition.
\begin{figure}[htbp]
	\centering
	\subfigure[$g=0.1$]{
		\includegraphics[width=4cm]{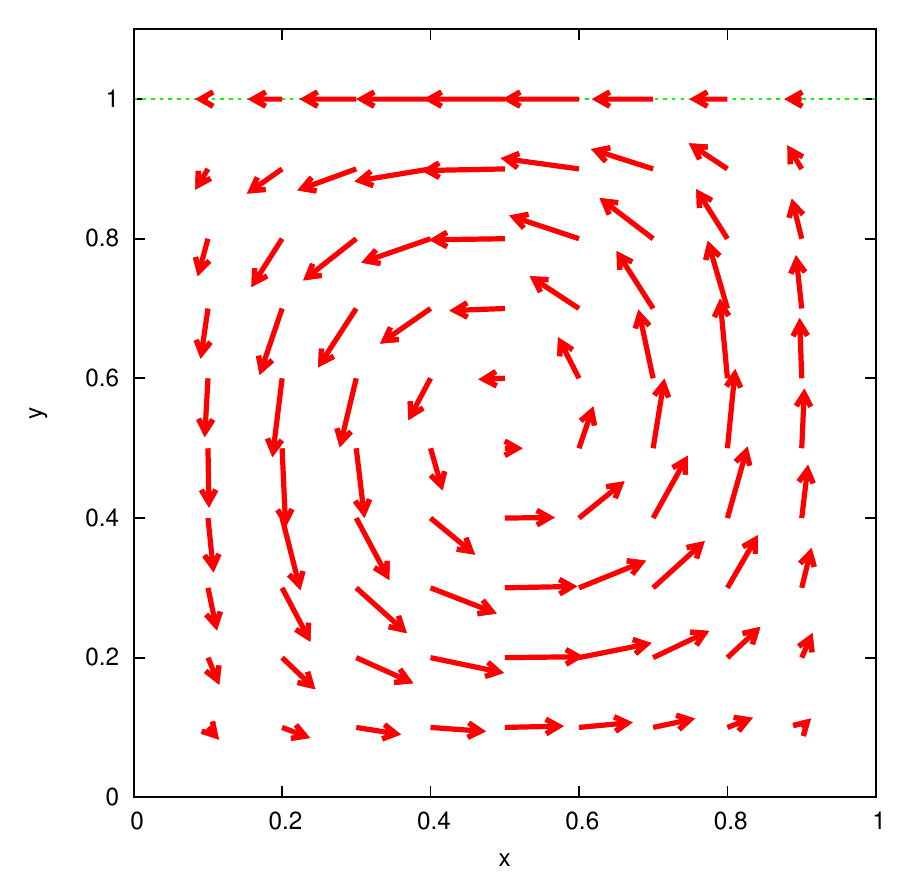}
	}
	\hspace{-5mm}
	\subfigure[$g=0.8$]{
		\includegraphics[width=4cm]{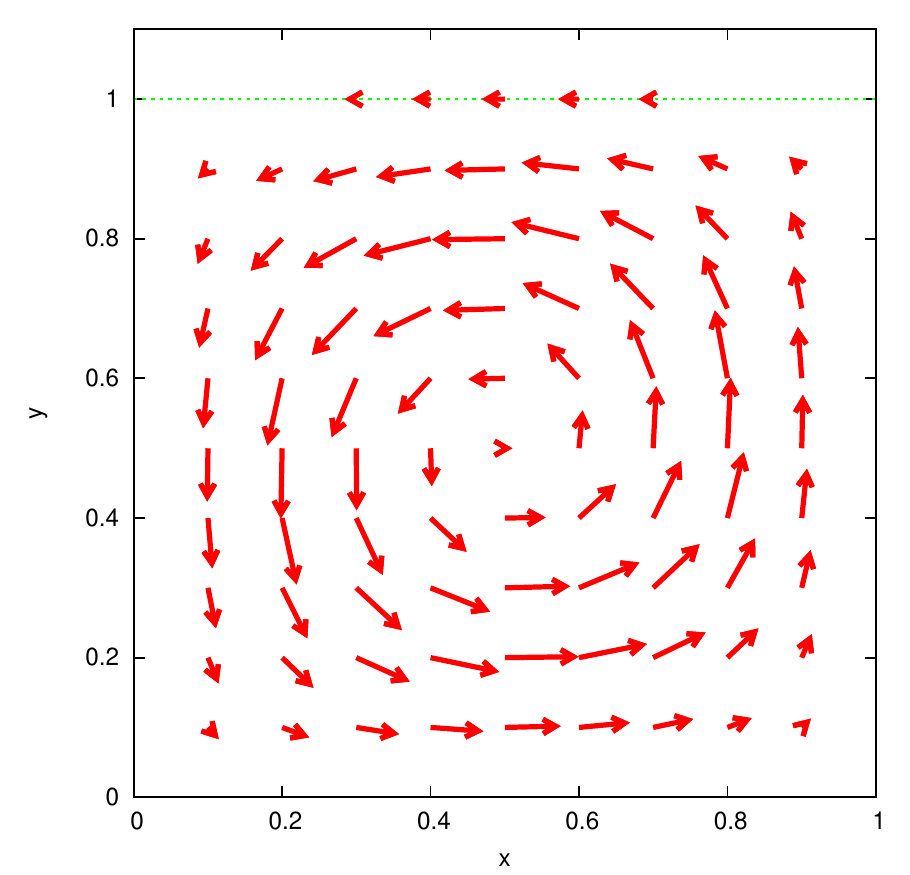}
	}
	\hspace{-5mm}
	\subfigure[$g=2.0$]{
		\includegraphics[width=4cm]{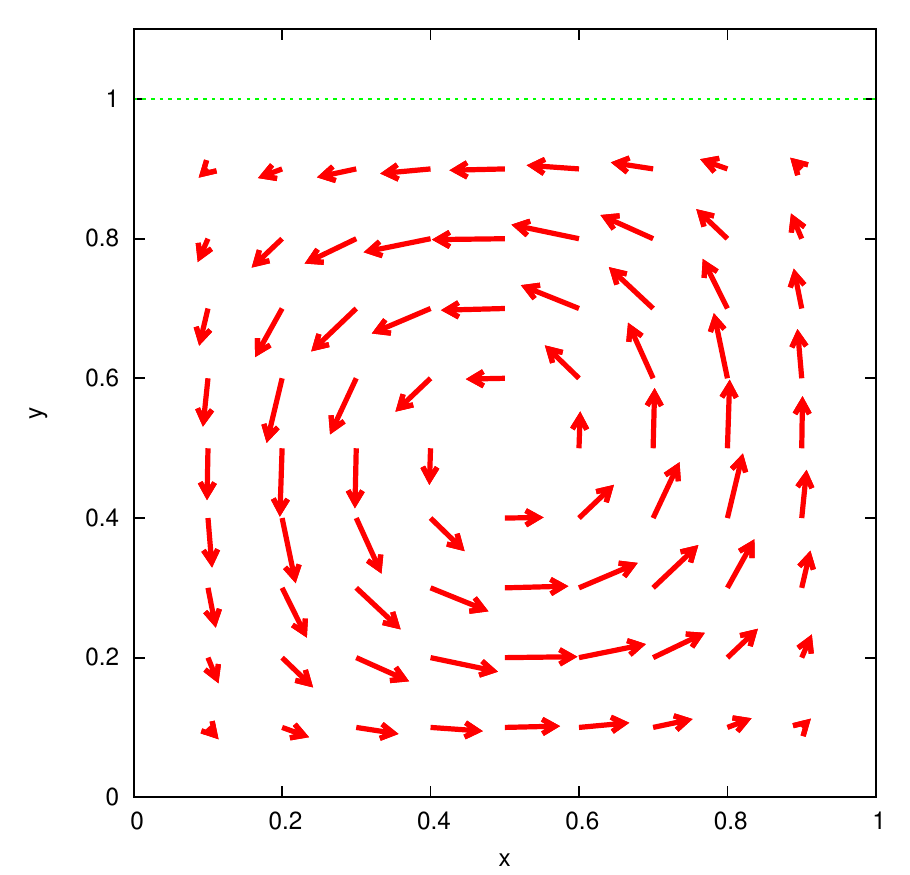}	
	}
	\caption{Solution velocity field of the Stokes equations with SBCF}
	\label{Fig5.1}
\end{figure}
\begin{figure}[htbp]
	\centering
	\subfigure[$g=0.1$]{
		\includegraphics[width=4cm]{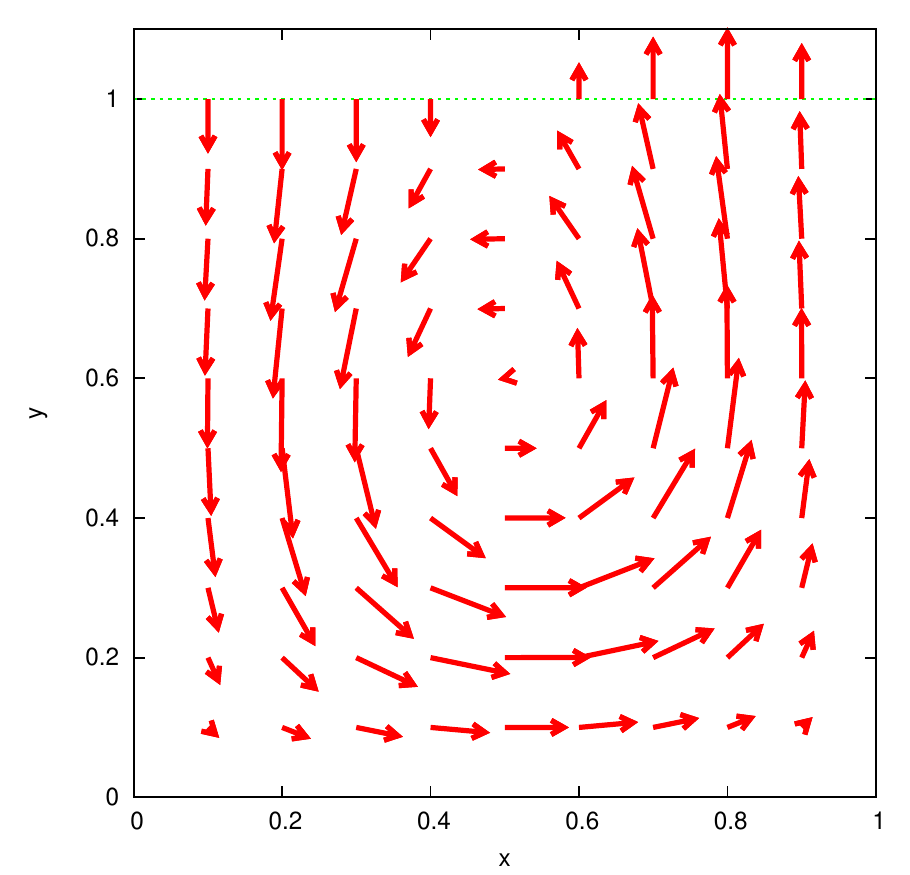}
	}
	\hspace{-5mm}
	\subfigure[$g=1.2$]{
		\includegraphics[width=4cm]{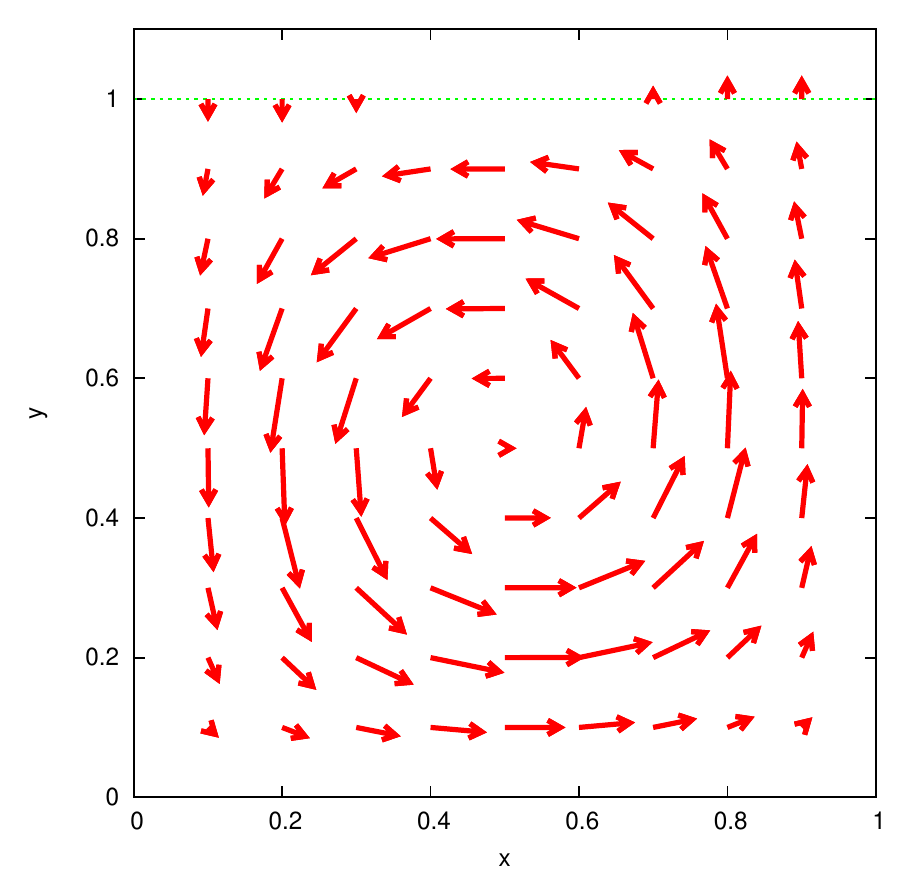}
	}
	\hspace{-5mm}
	\subfigure[$g=3.0$]{
		\includegraphics[width=4cm]{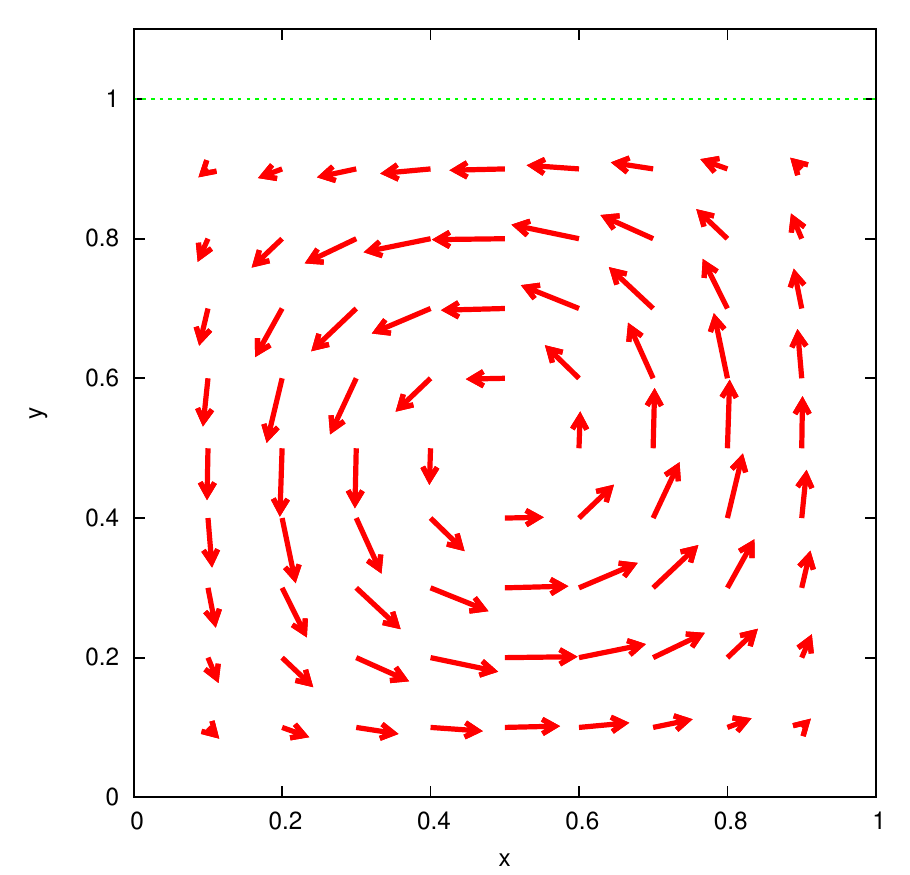}	
	}
	\caption{Solution velocity field of the Stokes equations with LBCF}
	\label{Fig5.2}
\end{figure}

\begin{table}
	\centering
	\caption{Values of the Lagrange multiplier $\lambda_h$ on $\Gamma_1$}
	\label{Tab5.1}
	\begin{tabular}{|c|ccc|ccc|c|}
	\hline
		 & \multicolumn{3}{c|}{SBCF} & \multicolumn{4}{c|}{LBCF} \\
	\hline
		$g$ & 0.1 & 0.8 & 2.0 & 0.1 & 1.2 & 3.0 & 3.0 \\
		$\rho$ & 1000.0 & 50.0 & 3.0 & 20.0 & 30.0 & 2.0 & 2.0 \\
		$\lambda_h^{(1)}$ & 0.0 & 0.0 & 0.0 & 0.0 & 0.0 & 0.0 & 0.2 \\ 
	\hline\hline
		$x$ &  &  &  &  &  &  & \\
	\hline
		0.0 & 0.0 & 0.0 & 0.0 & 0.0 & 0.0 & 0.0 & 0.0 \\
		0.1 & $-1.0$ & $-0.26$ & $-0.09$ & $-1.0$ & $-1.0$ & $-0.63$ & $-0.43$ \\
		0.2 & $-1.0$ & $-0.90$ & $-0.25$ & $-1.0$ & $-1.0$ & $-0.57$ & $-0.37$ \\
		0.3 & $-1.0$ & $-1.0$ & $-0.42$ & $-1.0$ & $-1.0$ & $-0.45$ & $-0.25$ \\
		0.4 & $-1.0$ & $-1.0$ & $-0.55$ & $-1.0$ & $-0.83$ & $-0.25$ & $-0.05$ \\
		0.5 & $-1.0$ & $-1.0$ & $-0.60$ & $-0.06$ & $-0.06$ & $-0.02$ & $0.18$ \\
		0.6 & $-1.0$ & $-1.0$ & $-0.55$ & $1.0$ & $0.67$ & $0.22$ & $0.42$ \\
		0.7 & $-1.0$ & $-1.0$ & $-0.43$ & $1.0$ & $1.0$ & $0.43$ & $0.63$ \\
		0.8 & $-1.0$ & $-0.94$ & $-0.26$ & $1.0$ & $1.0$ & $0.58$ & $0.78$ \\
		0.9 & $-1.0$ & $-0.26$ & $-0.09$ & $1.0$ & $1.0$ & $0.66$ & $0.86$ \\
		1.0 & 0.0 & 0.0 & 0.0 & 0.0 & 0.0 & 0.0 & 0.0 \\
	\hline
		$k_{\text{itr}}$ & 4 & 18 & 29 & 21 & 12 & 29 & 30 \\
	\hline
	\end{tabular}
	\vspace{1cm}
\end{table}

Next, we consider the behavior of the Lagrange multiplier $\lambda_h$.
It follows from (\ref{3.12}) or (\ref{4.15}) together with Lemma \ref{Lem3.3}(ii) that for each $M\in\; \stackrel\circ\Gamma_{1,h}$
\begin{equation}
	\begin{cases}
		|\lambda_h(M)| \le 1 &\text{ if }u_h(M)=0, \\
		\lambda_h(M) = +1 \text{ or } -1 &\text{ if }u_h(M)\neq0,
	\end{cases}
\end{equation}
which is observed by comparing the result of Table \ref{Tab5.1} with Figure \ref{Fig5.1} or \ref{Fig5.2}.
In Table \ref{Tab5.1}, we also see that if any leak does not occur, then the choice of the starting value $\lambda_h^{(1)}$ affects the resulting limit
($\lambda_h^{(1)}=0.2$ in the last column implies that $\lambda_h^{(1)}(M) = 0.2$ for each $M\in\stackrel\circ \Gamma_{1,h}$), whereas changing the value of $\rho$ does not cause such phenomena.
Here, all the computations shown in Figures \ref{Tab5.1}--\ref{Tab5.2} and Table \ref{Tab5.1} are performed for $N=10$ until the stopping criterion
\begin{equation}
	\|u_h^{(k)} - u_h^{(k-1)}\|_{H^1(\Omega)^2} \le 10^{-5} \label{5.2}
\end{equation}
is satisfied in Algorithm \ref{Alg3.1} or \ref{Alg4.1}. The number of iterations required to attain (\ref{5.2}) is denoted by $k_{\text{itr}}$.

\begin{table}
	\centering
	\caption{Convergence behavior of $\|u_h - u_{\rm ref}\|_{H^1(\Omega)^2}$ and $\|p_h - p_{\rm ref}\|_{L^2(\Omega)}$}
	\label{Tab5.2}
	\begin{tabular}{|c|cc|cc|cc|cc|}
	\hline
		$N$ & \multicolumn{4}{c|}{SBCF} & \multicolumn{4}{c|}{LBCF} \\
	\hline
		& $H^1$-error & rate & $L^2$-error & rate & $H^1$-error & rate & $L^2$-error & rate \\
	\hline\hline
		10 & 1.6E-2 & --- & 1.6E-2 & --- & 1.4E-2 & --- & 1.3E-2 & --- \\
		12 & 1.1E-2 & 1.9 & 1.1E-2 & 2.0 & 1.0E-2 & 1.8 & 9.7E-3 & 1.8 \\
		15 & 7.0E-3 & 2.1 & 6.3E-3 & 2.5 & 6.4E-3 & 2.0 & 5.8E-3 & 2.2 \\
		20 & 3.9E-3 & 2.0 & 3.5E-3 & 2.1 & 3.7E-3 & 1.9 & 3.3E-3 & 1.9 \\
		24 & 2.6E-3 & 2.1 & 2.7E-3 & 1.3 & 2.5E-3 & 2.2 & 2.2E-3 & 2.2 \\
		30 & 1.7E-3 & 2.0 & 1.5E-3 & 2.6 & 1.6E-3 & 2.0 & 1.5E-3 & 1.9 \\
		40 & 9.0E-4 & 2.1 & 8.5E-4 & 2.0 & 8.4E-4 & 2.2 & 8.0E-4 & 2.2 \\
	\hline
	\end{tabular}
\end{table}

Finally, we evaluate the error between approximate solutions and exact ones as the division number $N$ increases, when $g=0.8$ and $g=1.2$ for the case of SBCF and LBCF, respectively.
Since we do not know the explicit exact solutions, we employ the approximate solutions with $N=120$ as the reference solutions $(u_{\rm ref}, p_{\rm ref})$,
and numerically calculate $\|u_h - u_{\rm ref}\|_{H^1(\Omega)^2}$ and $\|p_h - p_{\rm ref}\|_{L^2(\Omega)}$.
Here, the additive constants of $p_h$'s are chosen such that $p_h(0,0) = p_{\text{ref}}(0,0)$.
Then, as Table \ref{Tab5.2} shows, we can observe the optimal order convergence $O(h^2)$ for both SBCF and LBCF.

\section{Conclusion and future works}\label{Sec7}
A finite element analysis using the P2/P1 element to the Stokes equations under SBCF or LBCF is examined.
We have proved the existence and uniqueness (partial non-uniqueness) results and established the convergence order $O(h^{1/4})$ as error estimates for appropriately smooth solutions;
sufficient conditions to obtain the optimal order $O(h^2)$ are also presented.
To compute the approximate solution, we have proposed an iterative Uzawa-type algorithm.
We have applied it to some examples and numerically observed the convergence order of $O(h^2)$.

In a future study, we would like to extend our theory to a more general situation, for example, a smooth domain without corners, nonlinear Navier-Stokes equations, a case in which SBCF and LBCF are imposed simultaneously, or a time-dependent problem.

\section*{Acknowledgments}
I would like to thank Dr. Hirofumi Notsu for providing a crucial idea that helped improve the performance of the numerical experiment.
I would also like to thank Professors Norikazu Saito and Hiroshi Suito for bringing this topic to my attention and encouraging me through valuable discussions.
This work was supported by CREST, JST.

\bibliographystyle{amsplain}
\providecommand{\bysame}{\leavevmode\hbox to3em{\hrulefill}\thinspace}
\providecommand{\MR}{\relax\ifhmode\unskip\space\fi MR }
\providecommand{\MRhref}[2]{%
  \href{http://www.ams.org/mathscinet-getitem?mr=#1}{#2}
}
\providecommand{\href}[2]{#2}

\end{document}